\documentclass[12pt,a4paper]{amsart}
\usepackage{amssymb}
\usepackage{latexsym}
\usepackage{amsmath}
\usepackage{amscd}
\usepackage[all]{xy}

\theoremstyle{plain}
\newtheorem{theorem}{Theorem}[section]
\newtheorem{prop}[theorem]{Proposition}
\newtheorem{cor}[theorem]{Corollary}
\newtheorem{lemma}[theorem]{Lemma}
\newtheorem{conj}[theorem]{Conjecture}
\theoremstyle{definition}

\newtheorem{defi}[theorem]{Definition}

\newtheorem{rem}[theorem]{Remark}
\newtheorem{remarks}[theorem]{Remarks}

\newcommand{\eqdef}{\;{:=}\;}

\newcommand{\ca}{\mathcal A}

\newcommand{\B}{{\mathbb B}}
\newcommand{\NN}{{\mathbb N}}
\newcommand{\C}{{\mathbb C}}
\newcommand{\Q}{{\mathbb Q}}
\newcommand{\R}{{\mathbb R}}

\newcommand{\Z}{{\mathbb Z}}

\newcommand{\V}{{\mathbb V}}

\newcommand{\p}{{\mathbb P}}

\newcommand{\op}{\operatorname}
\newcommand{\mc}[1]{{\mathcal #1}}

\newcommand{\SU}{\op{SU}}
\newcommand{\GU}{\op{GU}}
\newcommand{\SL}{\op{SL}}
\newcommand{\GL}{\op{GL}}

\newcommand{\End}{\op{End}}
\newcommand{\Aut}{\op{Aut}}
\newcommand{\Hom}{\op{Hom}}

\newcommand{\chow}[2]{{\rm CH}^{#1}(#2)}
\newcommand{\qchow}[2]{{\rm CH}^{#1}(#2)_{{\mathbb Q}}}
\newcommand{\corr}[4]{{\rm Corr}^{#1}(#3 \times_{#2} #4)}

\begin{document}
\parindent=0pt

\title[Motives for Picard modular surfaces]{Chow-K\"{u}nneth decomposition for universal families over Picard modular surfaces}
\author{ A. Miller, S. M\"{u}ller-Stach, S. Wortmann, Y.-H.Yang, K. Zuo}
\address{Mathematisches Institut der Universit\"{a}t Heidelberg, Im Neuenheimer Feld 288, 69120 Heidelberg }
\email{miller@mathi.uni-heidelberg.de, wortmann@mathi.uni-heidelberg.de}
\address{Mathematisches Institut der Johannes Gutenberg Universit\"at Mainz, Staudingerweg 9, 55099 Mainz}
\email{mueller-stach@mathematik.uni-mainz.de, kzuo@mathematik.uni-mainz.de}
\address{Max Planck Institut f\"ur Mathematik, Inselstrasse 22, 04103 Leipzig, Germany}
\email{yhyang@mail.tongji.edu.cn, yi-hu.yang@mis.mpg.de}

\thanks{Supported by: DFG Schwerpunkt--Programm, DFG China Exchange program, NSF of China (grant no. 10471105),
Max--Planck Gesellschaft.}
\date{\today}
\dedicatory{Dedicated to Jaap Murre}

\begin{abstract} We prove existence results for Chow--K\"unneth projectors on
compactified universal families of Abelian threefolds with complex
multiplication over a particular Picard modular surface studied by
Holzapfel. Our method builds up on the approach of Gordon, Hanamura
and Murre in the case of Hilbert modular varieties. In addition we
use relatively complete models in the sense of Mumford, Faltings and
Chai and prove vanishing results for $L^2$--Higgs cohomology groups
of certain arithmetic subgroups in $SU(2,1)$ which are not
cocompact.
\end{abstract}

\maketitle

\section{Introduction}

In this paper we discuss conditions for the existence of absolute
Chow-K\"unneth decompositions for families over Picard modular
surfaces and prove some partial existence results. In this way we
show how the methods of Gordon, Hanamura and
Murre~\cite{gordon-hanamura-murre-iii} can be slightly extended to
some cases but fail in some other interesting cases. Let us first
introduce the circle of ideas which are behind Chow--K\"unneth
decompositions. For a general reference we would like to encourage
the reader to look into
\cite{murre-motive} which gives a beautiful introduction to the subject and explains all notions we are using. \\
Let $Y$ be a smooth, projective $k$ --variety of dimension $d$ and
$H^*$ a Weil cohomology theory. In this paper we will mainly be
concerned with the case $k=\C$, where we choose singular cohomology
with rational coefficients as Weil cohomology. Grothendieck's
Standard Conjecture $C$ asserts that the K\"unneth components of the
diagonal $\Delta \subset Y \times Y$  in the cohomology $H^{2d}(Y
\times Y,\Q)$  are algebraic, i.e., cohomology classes of algebraic
cycles. In the case $k=\C$ this follows from the Hodge conjecture.
Since $\Delta$ is an element in the ring of correspondences, it is
natural to ask whether these algebraic classes come from algebraic
cycles $\pi_j$ which form a complete set of orthogonal idempotents
$$
\Delta = \pi_0 +\pi_1 + \ldots +\pi_{2d} \in CH^d(Y \times Y)_\Q
$$
summing up to $\Delta$. Such a decomposition is called a {\sl
Chow--K\"unneth decomposition} and it is conjectured to exist for
every smooth, projective variety. One may view $\pi_j$ as a Chow
motive representing the projection onto the $j$--the cohomology
group in a universal way. There is also a corresponding notion for
$k$--varieties which are relatively smooth over a base scheme $S$.
See section 3, where also Murre's refinement of this conjecture with
regard to the Bloch--Beilinson filtration is discussed.
Chow--K\"unneth decompositions for abelian varieties were first
constructed by Shermenev in 1974. Fourier--Mukai transforms may be
effectively used to write down the projectors, see
\cite{motives-kuenne, murre-motive}. The cases of surfaces was
treated by Murre~\cite{murre-surface}, in particular he gave a
general method to construct the projectors $\pi_1$ and $\pi_{2d-1}$,
the so--called Picard and Albanese Motives. Aside from other special
classes of $3$--folds \cite{dAMS} not much evidence is known except
for some classes of modular varieties. A fairly general method was
introduced and exploited recently by Gordon, Hanamura and Murre, see
\cite{gordon-hanamura-murre-iii}, building up on previous work by
Scholl and their own. It can be applied in the case where one has a
modular parameter space $X$ together with a universal family $f:A
\to X$ of abelian varieties with possibly some additional structure.
Examples are given by elliptic and Hilbert modular varieties. The
goal of this paper was to extend the range of examples to the case
of Picard modular surfaces, which are uniformized by a ball, instead
of a product of upper half planes. Let us now describe the general
strategy of Gordon, Hanamura and Murre so that we can understand to
what extent this approach differs and eventually fails
for a general Picard modular surface with sufficiently high level structure. \\
Let us assume that we have a family $f:A \to X$ of abelian varieties over $X$. Since all fibers are abelian,
we obtain a relative Chow---K\"unneth decomposition over $X$ in the sense of Deninger/Murre \cite{den-mur},
i.e., algebraic cycles $\Pi_j$ in
$A \times_X A$ which sum up to $\Delta \times_X \Delta$. One may view $\Pi_j$ as a projector related to
$R^jf_*\C$. Now let $\overline{f}: \overline{A} \to \overline{X}$ be a compactification of the family.
We will use the language of perverse sheaves from \cite{BBD} in particular also the notion of a stratified map.
In \cite{gordon-hanamura-murre-ii} Gordon, Hanamura and Murre have
introduced the {\sl Motivic Decomposition Conjecture} :

\begin{conj} \label{MDC} Let $\overline{A}$ and $\overline{X}$ be quasi--projective varieties over $\C$, $\overline{A}$
smooth, and $\overline{f}: \overline{A} \to \overline{X}$ a projective map. Let $\overline{X}=X_0 \supset X_1 \supset
\ldots \supset X_{\dim(X)}$ be a stratification of $\overline{X}$ so that $\overline{f}$ is a stratified map. Then there are local
systems ${\mathcal V}^j_\alpha$ on $X_\alpha^0=X_\alpha \setminus X_{\alpha-1}$, a complete set
$\Pi^j_\alpha$ of orthogonal projectors and isomorphisms
$$
\sum_{j,\alpha} \Psi^j_\alpha: \R\overline{f}_*\Q_{\overline{A}} {\buildrel \cong \over \to}
\bigoplus_{j,\alpha} IC_{X_\alpha}( {\mathcal V}^j_\alpha)[-j-\dim(X_\alpha)]
$$
in the derived category.
\end{conj}

This conjecture asserts of course more than a relative Chow--K\"unneth decomposition for the smooth part $f$
of the morphism $\overline{f}$. Due to the complicated structure of the strata in general its proof in general needs some more
information about the geometry of the stratified morphism $\overline {f}$.
In the course of their proof of the Chow--K\"unneth decomposition for Hilbert modular varieties, see \cite{gordon-hanamura-murre-iii},
Gordon, Hanamura and Murre have proved the motivic decomposition conjecture in the case of toroidal compactifications for the
corresponding universal families. However to complete their argument they need the vanishing theorem of Matsushima--Shimura \cite{matsushima-shimura}.
This theorem together with the decomposition theorem \cite{BBD} implies that each relative projector $\Pi^j$ on the generic stratum $X_0$
only contributes to one cohomology group of $A$ and therefore, using further reasoning on boundary strata $X_\alpha$,
relative projectors for the family $f$ already induce absolute projectors. \\
The plan of this paper is to extend this method to the situation of Picard modular surfaces.
These were invented by Picard in his study of the family of curves (called Picard curves) with the affine equation
$$
y^3=x(x-1)(x-s)(x-t).
$$
The Jacobians of such curves of genus $3$ have some additional $CM$--structure arising from the $\Z/3\Z$ deck transformation group.
Picard modular surfaces are compactifications of two dimensional ball quotients $X=\B_2/\Gamma$ which parametrize such Jacobians 
and form a particular beautiful set of Shimura surfaces in the moduli space of abelian varieties of dimension $3$. 
Many examples are known through the work of Holzapfel~\cite{holzapfel1, holzapfel2}.
Unfortunately the generalization of the vanishing theorem of
Matsushima and Shimura does not hold for Picard modular surfaces and their compactifications. The reason is that $\B_2$ is a homogenous space for the Lie group
$G=SU(2,1)$
 and general vanishing theorems like Ragunathan's theorem \cite[pg. 225]{borel-wallach} do not hold.
If $\V$ is an irreducible, non--trivial representation of any arithmetic subgroup $\Gamma$ of $G$, then the intersection
cohomology group $H^1(X, \V)$ is frequently  non--zero, whereas in order to make
the method of Gordon, Hanamura and Murre work, we would need its vanishing. This happens frequently
for small $\Gamma$, i.e., high level. However if $\Gamma$ is sufficiently big, i.e., the level is small,
we can sometimes expect some vanishing theorems to hold. This is the main reason why we concentrate our investigations
on one particular example of a Picard modular surface $\overline{X}$ in section 4.
The necessary vanishing theorems are proved by using Higgs bundles and their $L^2$--cohomology in section 6. Such techniques
provide a new method to compute intersection cohomology in cases where the geometry is known. This methods uses a recent proof of the Simpson
correspondence in the non--compact case by Jost, Yang and Zuo~\cite[Thm. A/B]{jyz}. But even in the case of our chosen surface $\overline{X}$
we are not able to show the complete vanishing result which would be necessary to proceed with the argument of
Gordon, Hanamura and Murre. We are however able to prove the existence of a partial set
$\pi_0, \pi_1, \pi_2, \pi_3, \pi_7, \pi_8, \pi_9, \pi_{10}$
of orthogonal idempotents under the assumption of the motivic decomposition conjecture \ref{MDC} on the universal family
$\overline{A}$ over $\overline{X}$:

\begin{theorem} \label{MT}
Assume the motivic decomposition conjecture \ref{MDC} for $\overline{f}: \overline{A} \to \overline{X}$.
Then ${\overline \ca}$ supports a partial set of Chow--K\"unneth projectors $\pi_i$ for $i \neq 4,5,6$.
\end{theorem}

Unfortunately we cannot prove the existence of the projectors
$\pi_4$, $\pi_5$, $\pi_6$ due to the non--vanishing of a certain
$L^2$--cohomology group, in our case $H^1(X,S^2{\mathbb V}_1)$,
where ${\mathbb V}_1$ is (half of) the standard representation. This
is special to $SU(2,1)$ and therefore the proposed method has no
chance to go through for other examples involving ball quotients. If
$H^1(X,S^2{\mathbb V}_1)$ would vanish or consist out of algebraic
Hodge $(2,2)$--classes only, then we would obtain a complete
Chow--K\"unneth decomposition. This is an interesting open question
and follows from the Hodge conjecture, if all classes in
$H^1(X,S^2{\mathbb V}_1)$ would have Hodge type $(2,2)$. 
We also sketch how to prove the motivic decomposition conjecture in this particular
case, see section 7.2., however details will be published elsewhere.
This idea generalizes the method from
\cite{gordon-hanamura-murre-iii}, since the fibers over boundary
points are not anymore toric varieties, but toric bundles over
elliptic curves. We plan to publish the full details in a
forthcoming publication and prefer to
assume the motivic decomposition conjecture \ref{MDC} in this paper.  \\
The logical structure of this paper is as follows:\\
In section 2 we present notations, definitions and known results
concerning Picard Modular surfaces and the universal Abelian schemes
above them. Section 3 first gives a short introduction to Chow
Motives and the Murre Conjectures and then proceeds to our case in
paragraph 3.2.
The remainder of the paper will then be devoted to the proof of Theorem \ref{MT}: \\
In section 4 we give a description of toroidal degenerations of
families of Abelian threefolds with complex multiplication.  In
section 5 we describe the geometry of a class of Picard modular
surfaces which have been studied by Holzapfel. In section 6 we prove
vanishing results for intersection cohomology using the non--compact
Simpson type correspondence between the $L^2$--Higgs cohomology of
the underlying VHS and the $L^2$--de Rham cohomology resp.
intersection cohomology of local systems. In section 7 everything is
put together to prove the main theorem \ref{MT}. The appendix
(section 8) gives an explicit description of the $L^2$--Higgs
complexes needed for the vanishing results of section 6.

\section{The Picard modular surface}\label{sec:def}

In this section we are going to introduce the (non--compact) Picard modular
surfaces $X=X_\Gamma$ and the universal abelian scheme $\ca$ of
fibre dimension $3$ over $X.$ For proofs and further references we
refer to \cite{gordon_canonical-models}.

Let $E$ be an imaginary quadratic field with ring of integers
$\mathcal{O}_E.$ The Picard modular group is defined as follows.
Let $V$ be a $3$-dimensional $E$-vector space and $L\subset V$ be
an $\mathcal{O}_E$-lattice. Let $J:V\times V\to E$ be a
nondegenerate Hermitian form of signature $(2,1)$ which takes
values in $\mathcal{O}_E$ if it is restricted to $L\times L.$ Now
let $G'=\SU(J,V)/\mathbb{Q}$ be the special unitary group of
$(V,\phi).$ This is a semisimple algebraic group over $\Q$ and for
any $\Q$-algebra $A$ its group of $A$-rational points is
\[G'(A)=\{g\in\SL(V\otimes_\Q A)\mid J(gu,gv)=J(u,v),\,\text{for all}\, u,v\in V\otimes_\Q A \}. \]
In particular one has $G'(\R)\simeq\SU(2,1)$.
The symmetric domain $\mathcal{H}$ associated to $G'(\R)$ can be
identified with the complex $2$-ball as follows. Let us fix once and
for all an embedding $E \hookrightarrow \C$ and identify $E\otimes_\Q \R$ with
$\C.$ This gives $V(\R)$ the structure of a $3$-dimensional
$\C$-vector space and one may choose a basis of $V(\R)$ such that
the form $J$ is represented by the diagonal matrix $[1,1,-1].$ As
$\mathcal{H}$ can be identified with the (open) subset of the
Grassmannian $\rm{Gr}_1(V(\R))$ of complex lines on which $J$ is
negative definite, one has
\[\mathcal{H}\simeq\{(Z_1,Z_2,Z_3)\in\mathbb{C}^3\mid
|Z_1|^2+|Z_2|^2-|Z_3|^2 <0 \}/\mathbb{C}^*.\] This is contained in
the subspace, where $Z_3\neq 0$ and, switching to affine
coordinates, can be identified with the complex $2$-ball
\[ \B=\{(z_1,z_2) \in {\mathbb C}^2 \mid
|z_1|^2+|z_2|^2<1 \}.\] Using this description one sees that
$G'(\R)$ acts transitively on $\B.$

The Picard modular group of $E$ is defined to be
$G'(\Z)=\SU(J,L),$ i.e. the elements $g\in G'(\Q)$ with $gL=L.$ It
is an arithmetic subgroup of $G(\mathbb{R})$ and acts properly
discontinuously on $\B.$ The same holds for any commensurable
subgroup $\Gamma\subset G'(\Q),$ in particular if $\Gamma\subset
G'(\Z)$ is of finite index the quotient
$X_\Gamma(\C)=\B/\Gamma$ is a non-compact complex
surface, the Picard modular surface. Moreover, for torsionfree
$\Gamma$ it is smooth.

We want to describe $X_\Gamma(\C)$ as moduli space for polarized
abelian $3$-folds with additional structure. For this we will give
a description of $X_\Gamma(\C)$ as the identity component of the
Shimura variety $S_K(G,\mathcal{H}).$

Let $G=\GU(J,V)/\Q$ be the reductive algebraic group of unitary
similitudes of $J,$ i.e. for any $\Q$-algebra $A$
\begin{multline*}
G'(A)=\{g\in\GL(V\otimes_\Q A)\mid \text{there
exists}\,\mu(g)\in A^* \, \text{such that} \\
J(gu,gv)=\mu(g)J(u,v),\,\text{for all}\, u,v\in V\otimes_\Q A\}.
\end{multline*}

As usual $\mathbb{A}$ denotes the $\Q$-adeles and $\mathbb{A}_f$
denotes the finite adeles. Let $K$ be a compact open subgroup of
$G(\mathbb{A}_f),$ which is compatible with the integral structure
defined by the lattice $L.$ I.e., $K$ should be in addition a
subgroup of finite index in $G(\hat{\Z})\eqdef\{g\in
G(\mathbb{A}_f)\mid g(L\otimes_\Z \hat{\Z})=L\otimes_\Z
\hat{\Z}\}.$ Then one can define
\[S_K(G,\mathcal{H})(\C)\eqdef G(\Q)\backslash \mathcal{H}\times G(\mathbb{A}_f)/K.\]
This can be decomposed as
$S_K(G,\mathcal{H})(\C)=\coprod_{j=1}^{n(K)}X_{\Gamma_j}(\C).$

The variety $S_K(G,\mathcal{H})(\C)$ has an interpretation as a
moduli space for certain $3$-dimensional abelian varieties. Recall
that over $\C$ an abelian variety $A$ is determined by the following
datum: a real vector space $W(\R),$ a lattice $W(\Z)\subset W(\R),$
and a complex strucuture $j:\C^\times\to \Aut_\R(W(\R))),$ for which
there exists a nondegenerate $\R-$bilinear skew-symmetric form
$\psi:W(\R)\times W(\R)\to \R$ taking values in $\Z$ on $W(\Z)$ such
that the form given by $(w,w')\mapsto \psi(j(i)w,w')$ is symmetric
and positive definite. The form $\psi$ is called a Riemann form and
two forms $\psi_1,\psi_2$ are called equivalent if there exist
$n_1,n_2\in\NN_{>0}$ such that $n_1\psi_1=n_2\psi_2.$ An equivalence
class of Riemann forms is called a homogeneous polarization of $A.$

An endomorphism of a complex abelian variety is an element
of\linebreak[4] $\End_\R(W(\R))$ preserving $W(\Z)$ and commuting
with $j(z)$ for all $z\in\C^\times.$ A homogenously polarized
abelian variety $(W(\R),W(\Z),j,\psi)$ is said to have complex
multiplication by an order $\mathcal{O}$ of $E$ if and only if
there is a homomorphism $m:\mathcal{O}\to\End(A)$ such that
$m(1)=1,$ and which is compatible with $\psi,$ i.e.
$\psi(m(\alpha^\rho)w,w')=\psi(w,m(\alpha)w')$ where $\rho$ is the
Galois automorphism of $E$ induced by complex conjugation (via our
fixed embedding $E \hookrightarrow \C$.) We shall only consider the case
$\mathcal{O}=\mathcal{O}_E$ in the following.

One can define the signature of the complex multiplication $m$,
resp. the abelian variety $(W(\R),W(\Z),j,\psi,m)$ as the
signature of the hermitian form $(w,w')\mapsto
\psi(w,iw')+i\psi(w,w')$ on $W(\R)$ with respect to the complex
structure imposed by $m$ via $\mathcal{O}\otimes_\Z \R\simeq \C.$
We write $m_{(s,t)}$ if $m$ has signature $(s,t).$

Finally for any compact open subgroup $K\subset G(\hat{\Z})$ as
before one has the notion of a level-$K$ structure on $A.$ For a
positive integer $n$ we denote by $A_n(\C)$ the group of points of
order $n$ in $A(\C).$ This group can be identified with
$W(\Z)\otimes\Z/n\Z$ and taking the projective limit over the system
$(A_n(\C))_{n\in\NN_{>0}}$ defines the Tate module of $A:$
\[T(A)\eqdef\varprojlim A_n(\C)\simeq W(\Z)\otimes\hat{\Z}.\]

Now two isomorphisms
$\varphi_1,\varphi_2:W(\Z)\otimes\hat{Z}\simeq L\otimes\hat{\Z}$
are called $K$-equivalent if there is a $k\in K$ such hat
$\varphi_1=k\varphi_2$ and a $K$-level structure on A is just a
$K$-equivalence class of these isomorphisms.

\begin{prop}
For any compact open subgroup $K\subset G(\hat{\Z})$ there is a
one-to-one correspondence between
\begin{enumerate}
    \item the set of points of $S_K(G,\mathcal{H})(\C)$ and
    \item the set of isomorphism classes of
    $(W(\R),W(\Z),j,\psi,m_{(2,1),\varphi})$ as above.
\end{enumerate}
\end{prop}
\begin{proof}
\cite[Prop.3.2]{gordon_canonical-models}
\end{proof}

\begin{rem}
If we take
$$
K_N\eqdef\{g\in G(\mathbb{A}_f)\mid (g-1)(L\otimes_\Z
\hat{\Z})\subset N\cdot(L\otimes_\Z \hat{\Z})\},
$$
then a level--$K$ structure is just the usual level-$N$ structure, namely an
isomorphism $A_N(\C)\to L\otimes \Z/N\Z.$ Moreover $K_N\subset
G(\Q)=\Gamma_N,$ where $\Gamma_N$ is the principal congruence
subgroup of level $N,$ i.e. the kernel of the canonical map
$G'(\Z)\to G'(\Z/N\Z).$ In this case the connected component of
the identity of $S_{K_N}(G,\mathcal{H})$ is exactly
$X_{\Gamma_N}(\C).$
\end{rem}

We denote with $\ca_\Gamma$ the universal abelian scheme over
$X_\Gamma(\C).$ In section \ref{sec:compact} the compactifications
of these varieties will be explained in detail. For the time being
we denote them with $\overline{X}_\Gamma$ and $\overline \ca _\Gamma.$

As the group $\Gamma$ will be fixed throughout the paper we will
drop the index $\Gamma$ if no confusion is possible.

\section{Chow motives and the conjectures of Murre}\label{sec:conj}
Let us briefly recall some definitions and results from the theory
of Chow motives. We refer to \cite{murre-motive} for details.
\subsection{}
For a smooth projective variety $Y$ over a field $k$ let
$\chow{j}{Y}$ denote the Chow group of algebraic cycles of
codimension $j$ on $Y$ modulo rational equivalence, and let
$\qchow{j}{Y}\eqdef\chow{j}{Y}\otimes\Q.$ For a cycle $Z$ on $Y$
we write $[Z]$ for its class in $\chow{j}{Y}.$ We will be working
with relative Chow motives as well, so let us fix a smooth
connected, quasi-projective base scheme $S\to \op{Spec}k.$ If
$S=\op{Spec}k$ we will usually omit $S$ in the notation. Let
$Y,Y'$ be smooth projective varieties over $S$, i.e., all fibers are smooth.
For ease of notation (and as we will not consider more general cases) we may
assume that $Y$ is irreducible and of relative dimension $g$ over
$S.$ The group of relative correspondences from $Y$ to $Y'$ of
degree $r$ is defined as
\[\corr{r}{S}{Y}{Y'}\eqdef \qchow{r+g}{Y\times_S Y'}.\]
Every $S$-morphism $Y'\to Y$ defines an element in
$\corr{0}{S}{Y}{Y'}$ via the class of the transpose of its graph.
In particular one has the class
$[\Delta_{Y/S}]\in\corr{0}{S}{Y}{Y}$ of the relative diagonal.
The self correspondences of degree $0$ form a ring, see \cite[pg. 127]{murre-motive}.
Using the relative correspondences one proceeds as usual to define
the category $\mathcal{M}_S$ of (pure) Chow motives over $S.$ The
objects of this pseudoabelian $\Q$-linear tensor category are
triples $(Y,p,n)$ where $Y$ is as above, $p$ is a projector, i.e.
an idempotent element in $\corr{0}{S}{Y}{Y},$ and $n\in\Z.$ The
morphisms are
\[\Hom_{\mathcal{M}_S}((Y,p,n),(Y',p',n'))
\eqdef p'\circ \corr{n'-n}{S}{Y}{Y'}\circ p.\]
When $n=0$ we write $(Y,p)$ instead of $(Y,p,0),$ and $h(Y)\eqdef
(Y,[\Delta_Y]).$

\begin{defi}
For a smooth projective variety $Y/k$ of dimension $d$ a
{\sl Chow-K{\"u}nneth-decomposition} of $Y$ consists of a collection of
pairwise orthogonal projectors $\pi_0,\ldots,\pi_{2d}$ in
$\corr{0}{}{Y}{Y}$ satisfying
\begin{enumerate}
\item $\pi_0+\ldots +\pi_{2d}=[\Delta_Y]$ and
\item for some Weil cohomology theory $H^*$ one has $\pi_i(H^*(Y))=H^i(Y).$
\end{enumerate}
\end{defi}
If one has a Chow-K{\"u}nneth decomposition for $Y$ one writes
$h^i(Y)=(Y,\pi_i)$. A similar notion of   a {\sl relative
Chow-K{\"u}nneth-decomposition} over $S$ can be defined in a straightforward manner, see also introduction.

Towards the existence of such decomposition one has the following
conjecture of Murre:

\begin{conj}
Let $Y$ be a smooth projective variety of dimension $d$ over some
field $k.$
\begin{enumerate}
\item There exists a Chow-K\"unneth decomposition for $Y$.
\item For all $i<j$ and $i>2j$ the action
of $\pi_i$ on $CH^j(Y)_\Q$ is trivial, i.e. $\pi_i\cdot CH^j(Y)_\Q=0$.
\item The induced $j$ step filtration on
$$
F^\nu CH^j(Y)_\Q:= {\rm Ker}\pi_{2j} \cap \cdots \cap
{\rm Ker}\pi_{2j-\nu+1}
$$
is independent of the choice of the
Chow--K\"unneth projectors, which are in general not canonical.
\item The first step of this filtration should give exactly the
subgroup of homological trivial cycles $CH^j(Y)_\Q$ in $CH^j(Y)_\Q$.
\end{enumerate}
\end{conj}

There are not many examples for which these conjectures have been
proved, but they are known to be true for surfaces
\cite{murre-motive}, in particular we know that we have a
Chow-K{\"u}nneth decomposition for $\overline{X}.$

In the following theorem we are assuming the motivic decomposition conjecture which was explained in the introduction.
The main result we are going to prove in section 7 is:

\begin{theorem}
Under the assumption of the motivic decomposition conjecture \ref{MDC}
$\overline \ca$ has a partial Chow--K\"unneth decomposition, including the projectors
$\pi_i$ for $i \neq 4,5,6$ as in Part (1) of Murre's conjecture.
\end{theorem}

Over the open smooth part $X\subset \overline{X}$ one has the relative projectors constructed
by Deninger and Murre in \cite{den-mur}, see also \cite{motives-kuenne}:
Let $S$ be a fixed base scheme as in section~\ref{sec:conj}. We will now
state some results on relative Chow motives in the case that $A$
is an abelian scheme of fibre dimension $g$ over $S.$
Firstly we have a functorial decomposition of the relative
diagonal $\Delta_{A/S}.$
\begin{theorem}
\label{uniquedec}
There is a unique decomposition
\[\Delta_{A/S}= \sum_{s=0}^{2g}\Pi_i\qquad\text{in}\qquad \qchow{g}{A\times_S A}\]
such that $(id_A\times [n])^*\Pi_i=n^i\Pi_i$ for all $n\in\Z.$
Moreover the $\Pi_i$ are mutually orthogonal idempotents, and
$[{}^t\Gamma_{[n]}]\circ\Pi_i=
n^i\Pi_i=\Pi_i\circ[{}^t\Gamma_{[n]}],$ where $[n]$ denotes the
multiplication by $n$ on $A.$
\end{theorem}
\begin{proof}
\cite[Thm.~3.1]{den-mur}
\end{proof}

Putting $h^i(A/S)=(A/S,\Pi_i)$ one has a Poincar\'{e}-duality for
these motives.

\begin{theorem}{\em (Poincar\'e-duality)}
\[h^{2g-i}(A/S)^\vee\simeq h^{i}(A/S)(g)\]
\end{theorem}
\begin{proof}
\cite[3.1.2]{motives-kuenne}
\end{proof}

\subsection{}

We now turn back to our specific situation. From Theorem
\ref{uniquedec} we have the decomposition
$\Delta_{\ca/X}=\Pi_0+\ldots +\Pi_6.$

We will have to extend these relative projectors to absolute
projectors. In order to show the readers which of the methods of
\cite{gordon-hanamura-murre-ii}, where Hilbert modular varieties are
considered, go through and which of them fail in our case, we recall
the main theorem (Theorem 1.3) from \cite{gordon-hanamura-murre-ii}:

\begin{theorem}
\label{ghm2} Let $p:\mathcal{A}\to X$ as above satisfy the following
conditions:

\begin{enumerate}
\item The irreducible components of $\overline{X} \setminus X$ are
smooth toric projective varieties.
\item The irreducible components of $\overline{\mathcal{A}} \setminus \mathcal{A}$
are smooth projective toric varieties.
\item The variety $\mathcal{A}/X$ has a relative Chow-K\"{u}nneth decomposition.
\item $\overline{X}$ has a Chow-K\"{u}nneth decomposition over $k$.
\item If $x$ is a point of $X$ the natural map.
\[CH^r(\mathcal{A})\to H^{2r}_B(\mathcal{A}_x(\mathbb{C}),\mathbb{Q})^{\pi_1^{top}(X,x)}\]
is surjective for $0\le r\le d=\dim \mathcal{A}-\dim X$.
\item For $i$ odd, $H^i_B(\mathcal{A}_x(\mathbb{C}),\mathbb{Q})^{\pi_1^{top}(X,x)}=0$.
\item  Let $\rho$ be an irreducible, non-constant
representation of ${\pi_1^{top}(X,x)}$ and $\mathbb{V}$ the
corresponding local system on $X$. Assume that $\mathbb{V}$ is
contained in the $i$--th exterior power $R^ip_*\Q=\Lambda^i
R^1p_*\Q$ of the monodromy representation for some $0 \le i \le 2d$.
Then the intersection cohomology $H^q(X,\mathbb{V})$ vanishes if
$q\neq\dim X$.
\end{enumerate}
Under these assumptions $\mathcal{A}$ has a Chow-K\"{u}nneth
decomposition over $k$.
\end{theorem}

As it stands we can only use conditions (3),(4) and (5) of this
theorem, all the other conditions fail in our case. As for
conditions (1) and (2) we will have to weaken them to torus
fibrations over an elliptic curve. This will be done in section 4.

Condition (3) holds because of the work of Deninger and Murre
(\cite{den-mur}) on Chow-K\"unneth decompositions of Abelian
schemes.

Condition (4) holds in our case because of the existence of
Chow-K\"unneth projectors for surfaces (see \cite{murre-motive}).

In order to prove condition (5) and to replace conditions (6) and
(7) we will from section 5 on use a non-compact Simpson type
correspondence between the $L^2$-Higgs cohomology of the underlying
variation of Hodge structures and the $L^2$-de Rham cohomology
(respectively intersection cohomology) of local systems. This will
show the vanishing of some of the cohomology groups mentioned in (6)
of Theorem \ref{ghm2} and enable us to weaken condition (7).

\section{The universal abelian scheme and its
compactification} \label{sec:compact}

In this section we show that the two conditions (1) and (2) of
Theorem \ref{ghm2} fail in our case. Instead of tori we get toric
fibrations over an elliptic curve as fibers over boundary
components. The main reference for this section is \cite{m2}.

\subsection{Toroidal compactifications of locally symmetric varieties}

In this paragraph an introduction to the theory of toroidal
compactifications of locally symmetric varieties as developed by
Ash, Mumford, Rapoport and Tai in \cite{amrt} is given. The main
goal is to fix notation. All details can be found in \cite{amrt},
see the page references in this paragraph.

\bigskip

Let $D=G(\R)/K$ be a bounded symmetric domain (or a finite number of
bounded symmetric domains, for the following discussion we will
assume $D$ to be just one bounded symmetric domain), where $G(\R)$
denotes the $\R$-valued points of a semisimple group $G$ and
$K\subset G(\R)$ is a maximal compact subgroup.  Let $\check{D}$ be
its compact dual. Then there is an embedding
\begin{equation}
\label{eqn:embedding} D\hookrightarrow\check{D}.
\end{equation}
Note that $G^0 (\C)$ acts on $\check{D}$.\\
We pick a parabolic $P$ corresponding to a rational boundary
component $F$, by $Z^0$ we denote the connected component of the
centralizer $Z(F)$ of $F$,
 by $P^0$ the connected component of $P$ and by $\Gamma$ a (torsion free, see below for this restriction) congruence subgroup of $G$. We will explicitly be interested only in connected groups, so from now on we can assume that $G^0=G$.  
Set
\[
\begin{split}
&\mbox{$N\subset P^0 $\quad  \quad  the unipotent radical}\\
&\mbox{$U\subset N $ \quad \quad the center of the unipotent radical}\\
&\mbox{$U_\C$ \quad  \quad  \quad  \quad    its complexification}\\
&V=N/U \\
&\Gamma_0=\Gamma\cap U\\
&\Gamma_1=\Gamma\cap P^0\\
&T=\Gamma_0\backslash U_\C.
\end{split}
\]
Note that $U$ is a real vector space and by construction, $T$ is an
algebraic torus over $\C$. Set $$D(F):=U_\C\cdot D,$$ where the dot
$\cdot$ denotes the action of  $G^0 (\C)$ on $\check{D}$. This is an
open set in $\check{D}$ and we have the inclusions
\begin{equation}
\label{eqn:inclusions} D\subset D(F)=U_\C\cdot D\subset \check{D}
\end{equation}
and furthermore a complex analytic isomorphism

\begin{equation}
\label{b} U_\C\cdot D\simeq U_\C\times\mc{E}_P
\end{equation}
where $\mc{E}_P$ is some complex vector bundle over the boundary
component corresponding to $P$. We will not describe $\mc{E}_P$ any
further, the interested reader is referred to \cite{amrt}, chapter
3. The isomorphism in \eqref{b} is complex analytic and takes the
$U_\C$-action on $\check{D}$ from the left to the translation on
$U_\C$ on the right.\\
We  once for all choose a boundary component $F$  and denote its
stabilizer by $P$.

\bigskip

From \eqref{eqn:inclusions} we get (see \cite{amrt}, chapter 3 for
details, e.g. on the last isomorphism)
\[
\Gamma_0\backslash D\subset \Gamma_0\backslash (U_\C\cdot D)\simeq
\Gamma_0\backslash(U_\C\times\mc{E}_P)\simeq T\times \mc{E},
\]
where $\mc{E}=\Gamma_0\backslash \mc{E}_P$.

The torus $T$ is the one we use for a toroidal embedding.
Furthermore $D$ can be realized as a Siegel domain of the third
kind:
\[
D\simeq \{(z,e)\in U_\C\cdot D\simeq U_\C\times\mc{E}\mid
\op{Im}(z)\in C+h(e)\},
\]
where
\[
h:\mc{E}\to U
\]
is a real analytic map and $C\subset U$ is an open cone in $U$.

A finer description of $C$ which is needed for the most general case
can be found in
 \cite{amrt}.

\bigskip

We pick a cone decomposition $\{\sigma_\alpha\}$ of $C$ such that

\begin{equation}
\begin{split}
\label{1/0}
&(\Gamma_1/\Gamma_0)\cdot \{\sigma_\alpha\}=\{\sigma_\alpha\} \text{ with finitely many orbits and }\\
&C \subset \bigcup_{\alpha}{\sigma_\alpha} \subset \overline{C}.
\end{split}
\end{equation}

This yields a torus embedding

\begin{equation}
\label{tembedding} T\subset X_{\{\sigma_\alpha\}}.
\end{equation}
We can thus partially compactify the open set $\Gamma_0\backslash
(U_\C\cdot D)$:
\[
\Gamma_0\backslash (U_\C\cdot D)\simeq\Gamma_0\backslash
U_\C\times\mc{E} \hookrightarrow X_{\{\sigma_\alpha\}}\times \mc{E}.
\]

\bigskip

The situation is now the following:
\begin{equation}
\label{T}
\begin{split}
\Gamma_0\backslash (U_\C\cdot D)& \simeq T\times\mc{E} \hookrightarrow X_{\{\sigma_\alpha\}}\times\mc{E}\\
{\cup}\quad&\\
\quad \quad \Gamma_0\backslash D.&
\end{split}
\end{equation}

We proceed to give a description of the vector bundle $\mc{E}_P$ in
order to describe the toroidal compactification geometrically.

Again from \cite{amrt} (pp.233) we know that
\[
\begin{split}
&D\cong F\times C\times N \qquad \text{as real manifolds and}\\
&D(F)\cong F \times V \times U_\C,
\end{split}
\]
where $V=N/U$ is the abelian part of $N$.


Now set
\[ D(F)':= D(F) \text{ mod } U_\C .
\]
This yields the following fibration:
\[
\xymatrix
{D(F) \ar@/_2pc/[dd]_\pi \ar[d]^{\text{fibres $U_\C$}}_{\pi_1}\\
D(F)' \ar[d]^{\text{fibres $V$}}_{\pi_2}\\
F.}
\]

Taking the quotient by $\Gamma_0$ yields a quotient bundle
\[
\xymatrix{
\Gamma_0\backslash D(F) \ar[d]^{\text{fibres $T:=\Gamma_0\backslash U_\C$}}_{\overline{\pi_1}}\\
D(F)'}
\]

So, $T$ is  an algebraic torus group with maximal compact subtorus
$$T_{cp}:=\Gamma_0\backslash U.$$

Take the closure of $\Gamma_0\backslash D$ in
$X_{\{\sigma_\alpha\}}\times\mc{E}$ and denote by
$(\Gamma_0\backslash D)_{\{\sigma_\alpha\}}$ its interior.

\bigskip

Factor $D\to \Gamma\backslash D$ by
\[
D\to \Gamma_0\backslash D\to \Gamma_1\backslash D\to
\Gamma\backslash D.
\]

It is the following situation we aim at obtaining:

\begin{equation}
\label{diagram}
\begin{array}{ccccccc}
(\Gamma_0\backslash D)_{\{\sigma_\alpha\}} & \hookleftarrow & \Gamma_0\backslash D & \to & \Gamma_1\backslash D & \to & \Gamma\backslash D\\
\cup &&\cup&&\cup&&\\
(\Gamma_0\backslash D(c))_{\{\sigma_\alpha\}} & \hookleftarrow &
\Gamma_0\backslash D(c) & \to & \Gamma_1\backslash D(c) &
\hookrightarrow & \Gamma\backslash D.
\end{array}
\end{equation}

Here $D(c)$ is a neighborhood of our boundary component.  More
precisely for any compact subset $K$ of the boundary component and
any $c\in C$ define
\[
D(c,K)=\Gamma_1\cdot \left\{
\begin{array}{ccc}
&& \op{Im} z \in C+h(e)+c\\
\raisebox{1.5ex}[-1.5ex]{$(z,e)\in U_\C\times\mc{E}$}&
\raisebox{1.5ex}[-1.5ex]{$\bigg|$} &\mbox{and $e$ lies above $K$}
\end{array}\right\}.
\]
Then by reduction theory for $c$ large enough, $\Gamma$-equivalence
on $D(c,K)$ reduces to $\Gamma_1$-equivalence.

\bigskip

This means that we have an inclusion
\[
\xymatrix{
\Gamma_1\backslash D(c)\quad \ar@{^{(}->}[r] \ar[d]^{\simeq} &\quad \Gamma\backslash D \\
(\Gamma_1/\Gamma_0)\backslash (\Gamma_0\backslash D(c))&}
\]
where the quotient by $\Gamma_1/\Gamma_0$ is defined in the obvious
way.

\bigskip

Furthermore $\Gamma_0\backslash D\hookrightarrow(\Gamma_0\backslash
D)_{\{\sigma_\alpha\}}$ directly induces
\[
\Gamma_0\backslash D(c)\hookrightarrow(\Gamma_0\backslash
D(c))_{\{\sigma_\alpha\}}.
\]
Having chosen $\{\sigma_\alpha\}$ such that
$(\Gamma_1/\Gamma_0)\cdot \{\sigma_\alpha\}=\{\sigma_\alpha\}$, we
get
\[
\Gamma_1\backslash D(c)\hookrightarrow (\Gamma_1/\Gamma_0)\backslash
(\Gamma_0\backslash D(c))_{\{\sigma_\alpha\}}
\]
which yields the partial compactification and establishes the
diagram (\eqref{diagram}).

The following theorem is derived from the above.

\begin{theorem}
\label{tor.comp} With the above notation and for a cone
decomposition $\{\sigma_\alpha\}$ of $C$ satisfying the condition
\eqref{1/0}, the diagram
\begin{equation}
\xymatrix{
\Gamma_1\backslash D(c) \quad\ar@{^{(}->}[r] \ar@{^{(}->}[dr] &\quad \Gamma\backslash D
\\
&(\Gamma_1/\Gamma_0)\backslash {(\Gamma_0 \backslash D(c))}_{ \{ \sigma_{\alpha} \} } }
\end{equation}
yields a (smooth if $\{\sigma_\alpha\}$ is chosen appropriately)
partial compactification of $\Gamma\backslash D$ at $F$.
\bigskip
\hfill $\Box$
\end{theorem}

\subsection{Toroidal compactification of Picard modular surfaces}

We will now apply the results of the last paragraph to the case of
Picard modular surfaces and give a finer description of the fibres
at the boundary.

\begin{theorem}
\label{complex} For each boundary component of a Picard modular
surface the following holds. With the standard notations from
\cite{amrt} (see also the last paragraph and \cite{m2}\, for the
specific choices of $\Gamma_0,\Gamma_1$ etc.)

$$(\Gamma_1/\Gamma_0)\backslash(\Gamma_0\backslash D)$$ is isomorphic to a punctured disc bundle over a CM elliptic curve $A$.

A toroidal compactification
\[
(\Gamma_1/\Gamma_0)\backslash (\Gamma_0\backslash D)_{\{\sigma_\alpha\}}
\]
is obtained by closing the disc with a copy of $A$ (e.g. adding the zero section of the corresponding line bundle).
\end{theorem}


We now turn to the modification of condition (2). The notation we
use is as introduced in chapter 3 of \cite{fc}. Let $\tilde P$ be a
relatively complete model of an ample degeneration datum associated
to our moduli problem. As a general reference for degenerations see
\cite{mum}, see \cite{fc} for the notion of relatively complete
model and \cite{m2} for the ample degeneration datum we need here.
In \cite{m2} the following theorem is proved.

\begin{theorem}
\label{thm:fibres}

\begin{itemize}
\item[(i)]
The generic fibre of $\tilde P$ is given by a
fibre-bundle over a CM elliptic curve $E$,  whose fibres are
countably many irreducible components of the form $\mathbb P$, where $\mathbb P$ is a
$\mathbb P^1$-bundle over $\mathbb P^1$.
\item[(ii)]
The special fibre of $\tilde P$ is given by a
fibre-bundle over the CM elliptic curve $E$,  whose fibres consist
of countably many irreducible components of the form $\mathbb P^1\times\mathbb P^1$.
\end{itemize}
\end{theorem}

\begin{rem} In this paper we work with some very specific Picard modular
surfaces and thus the generality of Theorem \ref{thm:fibres} is not
needed. It will be needed though to extend our results to larger
families of Picard modular surfaces, see section 7.2..
\end{rem}

\section{Higgs bundles on Picard modular surfaces}

In this section we describe in detail the Picard modular surface of Holzapfel which is our main object.
We follow Holzapfel~\cite{holzapfel1, holzapfel2} very closely. In the remaining part of this section we explain the
formalism of Higgs bundles which we will need later.

\subsection{Holzapfel's surface}

We restrict our attention to the \emph{Picard modular surfaces} with compactification $\overline X$
and boundary divisor $D \subseteq \overline X$ which were
discussed by Picard \cite{picard}, Hirzebruch \cite{hirzebruch} and Holzapfel \cite{holzapfel1,holzapfel2}.
These surfaces are compactifications of ball quotients $X={\mathbb B}/\Gamma$ where
$\Gamma$ is a subgroup of $SU(2,1;{\mathcal O})$ with
${\mathcal O}={\mathbb Z} \oplus {\mathbb Z} \omega$, $\omega=\exp(2 \pi i/3)$, i.e.,
${\mathcal O}$ is the ring of Eisenstein numbers. In the case $\Gamma=SU(2,1;{\mathcal O})$, studied already by Picard,
the quotient ${\mathbb B}/\Gamma$ is ${\mathbb P}^2 \setminus$ $4$ points, an open set of which
is $U={\mathbb P}^2 \setminus \Delta$ and $\Delta$ is a configuration of $6$ lines
(not a normal crossing divisor). $U$ is a natural parameter space for a family of \emph{Picard curves}
$$
y^3=x(x-1)(x-s)(x-t)
$$
of genus $3$ branched over $5$ (ordered) points $0,1,s,t,\infty$ in ${\mathbb P}^1$.
The parameters $s,t$ are coordinates in the affine set $U$. If one looks at the subgroup
$$
\Gamma'= \Gamma \cap SL(3,{\mathbb C}),
$$
then $X={\mathbb B}/\Gamma'$ has a natural compactification $\overline X$ with a smooth boundary divisor
$D$ consisting of $4$ disjoint elliptic curves $E_0+E_1+E_2+E_3$, see \cite{holzapfel1,holzapfel2}.
This surface $\overline X$ is birational to a covering of ${\mathbb P}^2 \setminus \Delta$
and hence carries a family of curves over it. If we pass to yet another subgroup
$\Gamma'' \subset \Gamma$ of finite index, then we obtain a Picard modular surface
$$
\overline X=\widetilde{\widetilde{E \times E}}
$$
with boundary $D$ a union of $6$ elliptic curves which are the strict transforms of the following $6$ curves
$$
T_1,T_\omega,T_{\omega^2},E \times \{Q_0\}, E \times \{Q_1\}, E \times \{Q_2\}
$$
on $E \times E$ in the notation of \cite[page 257]{holzapfel1}.
This is the surface we will study in this paper.
The properties of the modular group $\Gamma''$ are described in
\cite[remark V.5]{holzapfel1}. In particular it acts freely on the ball. $\overline X$ is the blowup of $E \times E$
in the three points $(Q_0,Q_0)$ (the origin), $(Q_1,Q_1)$ and $(Q_2,Q_2)$ of triple intersection.
Note that $E$ has the equation $y^2z=x^3-z^3$.
On $E$ we have an action of $\omega$ via $(x:y:z) \mapsto (\omega x:y:z)$. $E$ maps to $\p^1$
using the projection
$$
p:E \to \p^1, \quad (x:y:z) \mapsto (y:z).
$$
This action has $3$ fixpoints
$Q_0=(0:1:0)$ (the origin), $Q_1=(0:i:1)$ and $Q_2=(0:-i:1)$ which are triple ramification points
of $p$. Therefore one has $3Q_0=3Q_1=3Q_2$ in $CH^1(E)$.

In order to proceed, we need to know something about the Picard group of $\overline X$.

\begin{lemma} In ${\rm NS}(E \times E)$ one has the relation
$$
T_1+T_{\omega}+T_{\omega^2}=3(0 \times E)+3(E \times 0).
$$
\end{lemma}

\begin{proof} Since $E$ has complex multiplication by $\Z[\omega]$,
the N\'eron--Severi group has rank $4$ and
divisors $T_1$, $T_{\omega}$, $0 \times E$ and $E \times 0$ form a basis of
${\rm NS}(E \times E)$. Using the intersection matrix of this basis, the claim
follows.
\end{proof}

The following statement is needed later:

\begin{lemma} The log--canonical divisor is divisible by three:
$$
K_{\overline{X}}+D=3L
$$
for some line bundle $L$.
\end{lemma}

\begin{proof}
If we denote by $\sigma: \overline{X} \to E \times E$ the blowup in the three points
$(Q_0,Q_0)$, $(Q_1,Q_1)$ and $(Q_2,Q_2)$, then we denote by
$Z=Z_1+Z_2+Z_3$ the union of all exceptional divisors. We get:
$$
\sigma^* T = D_1+Z, \; \sigma^* T_{\omega} = D_2+Z, \;  \sigma^* T_{\omega^2} = D_3+Z,
$$
and
$$
\sigma^* E \times Q_0=D_4+Z_1, \; \sigma^* E \times Q_1=D_5+Z_2, \; \sigma^* E \times Q_2=D_6+Z_3.
$$
Now look at the line bundle $K_{\overline{X}}+D$. Since
$$
K_{\overline{X}}+D=\sigma^* K_{E \times E} + Z + D=Z+D,
$$
we compute
$$
K_{\overline{X}}+D=\sum_{i=1}^6 D_i +\sum_{j=1}^3 Z_j.
$$
The first sum,
$$
D_1+D_2+D_3=-3Z+\sigma^*(T_1+T_{\omega}+T_{\omega^2})=-3Z+3 \sigma^*(0 \times E + E \times 0).
$$
is divisible by $3$. Using $3Q_0=3Q_1=3Q_2$, the rest can be computed in ${\rm NS}(\overline{X})$ as
$$
D_4+D_5+D_6+Z= \sigma^*(E \times 0 + E\times Q_1 + E \times Q_2)=3 \sigma^*(E \times 0).
$$
Therefore the class of $K_{\overline{X}}+D$ in ${\rm NS}(\overline{X})$ is given by
$$
K_{\overline{X}}+D= -3Z+ 3 \sigma^*(0 \times E) + 6 \sigma^*(E \times 0)
$$
and divisible by $3$. Since ${\rm Pic}^0(\overline{X})$ is a divisible group, $K_{\overline{X}}+D$
is divisible by $3$ in ${\rm Pic}(\overline{X})$ and we get a line bundle $L$ with
$K_{\overline{X}}+D=3L$ whose class in ${\rm NS}(\overline{X})$ is given by
$$
L=\sigma^*(0 \times E)-Z + 2\sigma^*(E \times 0).
$$
If we write
$$
\sigma^*(0 \times E)=D_0+Z_1,
$$
we obtain the equation
$$
L=D_0+D_5+D_6
$$
in ${\rm NS}(\overline{X})$. Note that $D_0$ intersects both $D_5$ and $D_6$ in one point.
All $D_i$, $i=1,\ldots,6$ have negative selfintersection and are disjoint.
\end{proof}

It is not difficult to see that $L$ is a nef and big line bundle since $\overline{X}$
has logarithmic Kodaira dimension $2$ \cite{holzapfel1}.
$L$ is trivial on all components of $D$ by the adjunction formula,
since they are smooth elliptic curves. \\
\ \\
The rest of this section is about the rank $6$ local system ${\mathbb V}=R^1p _* {\mathbb Z}$
on $X$. The following Lemma was known to Picard~\cite{picard}, he wrote down
$3 \times 3$ monodromy matrices with values in the Eisenstein numbers:

\begin{lemma} $\V$ is a direct sum of two local systems $\V=\V_1 \oplus \V_2$ of rank $3$.
The decomposition is defined over the Eisenstein numbers.
\end{lemma}
\begin{proof}
The cohomology $H^1(C)$ of any Picard curve $C$ has a natural
${\mathbb Z}/3{\mathbb Z}$ Galois action. Since the projective line has $H^1({\mathbb P}^1,{\mathbb Z})=0$,
the local system ${\mathbb V} \otimes \C$ decomposes into two $3$--dimensional local systems
$$
{\mathbb V}={\mathbb V}_1 \oplus {\mathbb V}_2
$$
which are conjugate to each other and defined over the Eisenstein numbers.
\end{proof}

Both local systems $\V_1, \V_2$ are irreducible and non--constant.

\subsection{Higgs bundles on Holzapfel's surface}

Now we will study the categorical correspondence between local systems and {\sl Higgs bundles}.
It turns out that it is sometimes easier to deal with one resp. the other.

\begin{defi} A Higgs bundle on a smooth variety $Y$ is a holomorphic vector bundle $E$ together with
a holomorphic map
$$
\vartheta: E \to E \otimes \Omega^1_Y
$$
which satisfies $\vartheta \wedge \vartheta=0$, i.e., an ${\rm End}(E)$ valued holomorphic $1$--form on $Y$.
\end{defi}

Each Higgs bundle induces a complex of vector bundles:
$$
E \to E \otimes \Omega^1_Y \to E \otimes \Omega^2_Y \to \ldots \to E \otimes \Omega^d_Y.
$$
{\sl Higgs cohomology} is the cohomology of this complex.
The {\sl Simpson} correspondence on a projective variety $Y$ gives an equivalence of categories between
polystable Higgs bundles with vanishing Chern classes and semisimple local systems $\V$ on $Y$ \cite[Sect. 8]{simpson}.
This correspondence is very difficult to describe in general and uses a deep existence theorem for harmonic metrics.
For quasi--projective $Y$ this may be generalized provided
that the appropriate harmonic metrics exist, which is still not known until today. There is however the case of VHS (Variations of Hodge structures) where the harmonic metric is the Hodge metric and is canonically given. For example if we have a smooth, projective family
$f: A \to X$ as in our example and $\V=R^mf_*\C$ is a direct image sheaf, then the corresponding Higgs bundle is
$$
E=\bigoplus_{p+q=m} E^{p,q}
$$
where $E^{p,q}$ is the $p$--the graded piece of the Hodge filtration $F^\bullet$ on ${\mathcal H}=\V \otimes {\mathcal O}_X$.
The Higgs operator $\vartheta$ is then given by the graded part of the Gau{\ss}--Manin connection, i.e., the cup product with the
Kodaira--Spencer class. In the non--compact case there is also a corresponding log--version for Higgs bundles, where $\Omega^1_Y$ is replaced by
$\Omega^1_Y(\log D)$ for some normal crossing divisor $D \subset Y$ and $E$ by the Deligne extension. Therefore we have to assume that
the monodromies around the divisors at infinity are unipotent and not only quasi--unipotent as in \cite[Sect. 2]{jyz}. This is the case in
Holzapfel's example, in fact above we have already checked that the log--canonical divisor $K_{\overline{X}}+D$
is divisible by three. We refer to \cite{simpson} and \cite{jyz} for the general theory.
In our case let $E=E^{1,0} \oplus E^{0,1}$ be the Higgs bundle corresponding to $\V_1$ with Higgs field
$$
\vartheta: E \to E \otimes \Omega^1_{\overline{X}}(\log D).
$$
This bundle is called the \emph{uniformizing bundle} in \cite[Sect. 9]{simpson}. \\
Let us return to Holzapfel's example.
We may assume that $E^{1,0}$ is $2$--dimensional and $E^{0,1}$ is $1$--dimensional, otherwise we
permute $\V_1$ and $\V_2$.
\begin{lemma}
$\vartheta: E^{1,0} \to E^{0,1} \otimes \Omega^1_{\overline X}(\log D)$
is an isomorphism.
\end{lemma}
\begin{proof} For the generic fiber this is true for rank reasons. At the boundary $D$ this is a local computation using the definition of
the Deligne extension. This has been shown in greater generality in \cite[Sect. 2-4]{jyz} (cf. also \cite[Sect. 4]{looijenga}),
therefore we do not give any more details.
\end{proof}

Let us summarize what we have shown for Holzapfel's surface $\overline{X}$:

\begin{cor}
$K_{\overline X}(D)$ is nef and big and there is a nef and big line bundle $L$ with
$$
L^{\otimes 3} \cong K_{\overline X}(D).
$$
The uniformizing bundle $E$ has components
$$
E^{1,0}=\Omega^1_{\overline X}(\log D) \otimes L^{-1}, \quad E^{0,1} =L^{-1}.
$$
The Higgs operator $\vartheta$ is the identity as a map $E^{1,0} \to E^{0,1} \otimes \Omega^1_{\overline X}(\log D)$ and it is trivial on $E^{0,1}$.
\end{cor}

\section{Vanishing of intersection cohomology}

Let $X$ be Holzapfel's surface from the previous section.
We now want to discuss the vanishing of intersection cohomology
$$
H^1(X,{\mathbb W})
$$
for irreducible, non--constant local systems ${\mathbb W} \subseteq R^ip_*\Q$.
Let ${\mathbb V}_1$ be as in the previous section.
Denote by $(E,\vartheta)$ the corresponding Higgs bundle with
$$
E= \left( \Omega^1_{\overline X}(\log D) \otimes L^{-1} \right) \oplus L^{-1}
$$
and Higgs field
$$
\vartheta: E \to E \otimes \Omega^1_{\overline X}(\log D).
$$

Our goal is to compute the intersection cohomology of $\V_1$.
We use the isomorphism between $L^2$-- and intersection cohomology for $\C$--VHS, a theorem of Cattani, Kaplan
and Schmid together with the isomorphism between $L^2$--cohomology and $L^2$--Higgs cohomology from \cite[Thm. A/B]{jyz}.
Therefore for computations of intersection cohomology we may use $L^2$--Higgs cohomology. We refer to
\cite{jyz} for a general introduction to all cohomology theories involved.

\begin{theorem} \label{firstvanishing}
The intersection cohomology $H^q(X,{\mathbb V}_1)$ vanishes for $q \ne 2$. 
By conjugation the same holds for ${\mathbb V}_2$. 
\end{theorem}

\begin{proof}
We need only show this for $q=1$, since ${\mathbb V}_1$ has no invariant sections, hence
$H^0(X,{\mathbb V}_1)=0$ and the other vanishings follow via duality
$$
H^q(X,{\mathbb V}_1) \cong H^{2\dim(X)-q}(X,{\mathbb V}_2)
$$
from the analogous statement for ${\mathbb V}_2$.
The following theorem provides the necessary technical tool.
\end{proof}

\begin{theorem}[{\cite[Thm. B]{jyz}}]
The intersection cohomology $H^q(X,{\mathbb V}_1)$ can be computed as the $q$--th hypercohomology
of the complex
$$
0 \to \Omega^0(E)_{(2)} \to \Omega^1(E)_{(2)} \to \Omega^2(E)_{(2)} \to 0
$$
on $\overline X$, where $E$ is as above. This is a subcomplex of
$$
E {\buildrel \vartheta \over \to} E \otimes \Omega^1_{\overline X}(\log D)
{\buildrel \vartheta \over \to} E \otimes \Omega^2_{\overline X}(\log D).
$$
In the case where $D$ is smooth, this is a proper subcomplex with the property
$$
\Omega^1(E)_{(2)} \subseteq  \Omega^1_{\overline X} \otimes E.
$$
\end{theorem}

\begin{proof}
This is a special case of the results in \cite{jyz}. The subcomplex is explicitly described
in section~\ref{appendix} of our paper.
\end{proof}

\begin{lemma} Let $E$ be as above with $L$ nef and big. Then the vanishing
$$
H^0(\Omega^1_{\overline X}(\log D) \otimes \Omega^1_{\overline X} \otimes L^{-1})=0
$$
implies the statement of theorem~\ref{firstvanishing}.
\end{lemma}

\begin{proof} We first compute the cohomology groups for the complex of vector bundles
and discuss the $L^2$--conditions later.
Any logarithmic Higgs bundle $E=\oplus E^{p,q}$ coming from a VHS has differential
$$
\vartheta: E^{p,q} \to E^{p-1,q+1} \otimes \Omega^1_{\overline X}(\log D).
$$
In our case $E=E^{1,0} \oplus E^{0,1}$ and the restriction of
$\vartheta$ to $E^{0,1}$ is zero. The differential
$$
\vartheta: E^{1,0} \to E^{0,1} \otimes \Omega^1_{\overline X}(\log D)
$$
is the identity. Therefore the complex
$$
(E^\bullet,\vartheta): E {\buildrel \vartheta \over \to} E \otimes \Omega^1_{\overline X}(\log D)
{\buildrel \vartheta \over \to} E \otimes \Omega^2_{\overline X}(\log D)
$$
looks like:
\begin{tiny}
$$
\begin{matrix}
&&&& \left(\Omega^1_{\overline X}(\log D) \otimes L^{-1}\right) & \oplus & L^{-1} \cr
&&&& \downarrow \cong  && \downarrow   \cr
&&\left(\Omega^1_{\overline X}(\log D)^{\otimes 2} \otimes L^{-1} \right)
& \oplus & \left(L^{-1} \otimes \Omega^1_{\overline X}(\log D) \right) & & 0  \cr
&&\downarrow   &&&&  \cr
\left(\Omega^1_{\overline X}(\log D) \otimes L^{-1} \otimes \Omega^2_{\overline X}(\log D) \right) &  \oplus &
\left( L^{-1} \otimes \Omega^2_{\overline X}(\log D)\right). & & & &
\end{matrix}
$$
\end{tiny}
Therefore it is quasi--isomorphic to a complex
$$
L^{-1} {\buildrel 0 \over \longrightarrow} S^2 \Omega^1_{\overline X}(\log D) \otimes L^{-1} {\buildrel 0 \over \longrightarrow}
\Omega^1_{\overline X}(\log D) \otimes \Omega^2_{\overline X}(\log D) \otimes L^{-1}
$$
with trivial differentials. As $L$ is nef and big, we have
$$
H^0(L^{-1})=H^1(L^{-1})=0.
$$
Hence we get
$$
{\mathbb H}^1({\overline X},(E^\bullet,\vartheta))
\cong H^0({\overline X}, S^2 \Omega^1_{\overline X}(\log D) \otimes L^{-1})
$$
and ${\mathbb H}^2({\overline X},(E^\bullet,\vartheta))$ is equal to
\begin{small}
$$
H^0({\overline X},K_{\overline X} \otimes L)^\vee \oplus
H^0({\overline X},\Omega^1_{\overline X}(\log D) \otimes \Omega^2_{\overline X}(\log D) \otimes L^{-1})
\oplus H^1({\overline X}, S^2 \Omega^1_{\overline X}(\log D) \otimes L^{-1}).
$$
\end{small}
If we now impose the $L^2$--conditions and use the complex
$\Omega^*_{(2)}(E)$ instead of $(E^\bullet,\vartheta)$, the resulting cohomology groups are subquotients
of the groups described above. Since
$$
\Omega^1(E)_{(2)} \subseteq  \Omega^1_{\overline X} \otimes E
$$
we conclude that the vanishing
$$
H^0({\overline X}, \Omega^1_{\overline X}(\log D) \otimes \Omega^1_{\overline X} \otimes L^{-1})=0
$$
is sufficient to deduce the vanishing of intersection cohomology.
\end{proof}

Now we verify the vanishing statement.

\begin{lemma} In the example above we have
$$
H^0(\Omega^1_{\overline X}(\log D) \otimes \Omega^1_{\overline X} \otimes L^{-1})=0.
$$
\end{lemma}

\begin{proof} Let $\sigma: \overline X \to E \times E$ be the blow up of the $3$ points of intersection.
Then one has an exact sequence
$$
0 \to \sigma^* \Omega^1_{E \times E} \to \Omega^1_{\overline X} \to i_* \Omega^1_Z \to 0,
$$
where $Z$ is the union of all (disjoint) exceptional divisors. Now $\Omega^1_{E \times E}$ is the trivial sheaf of rank $2$.
Therefore $\Omega^1_{\overline X}(\log D) \otimes \Omega^1_{\overline X} \otimes L^{-1}$
has as a subsheaf $2$ copies of $\Omega^1_{\overline X}(\log D) \otimes L^{-1}$.
The group
$$
H^0(\overline X,  \Omega^1_{\overline X}(\log D) \otimes L^{-1})
$$
is zero by the Bogomolov--Sommese vanishing theorem (see \cite[Cor. 6.9]{ev}),
since $L$ is nef and big. In order to prove the assertion it is hence sufficient to show that
$$
H^0(Z, \Omega^1_{\overline X}(\log D) \otimes  \Omega^1_Z \otimes L^{-1} )=0.
$$
But $Z$ is a disjoint union of ${\mathbb P}^1$'s. In our example we have
$K_{\overline X}(D) \otimes {\mathcal O}_Z \cong {\mathcal O}_Z(3)$ since $(L.Z)=1$
and therefore $\Omega^1_Z \otimes L^{-1} \cong {\mathcal O}_Z(-3)$.
Now we use in addition the conormal sequence
$$
0 \to N^*_Z \to \Omega^1_{\overline X}(\log D)|_Z  \to \Omega^1_{Z}(\log (D \cap Z)) \to 0.
$$
Note that $N^*_Z={\mathcal O}_Z(1)$.
Twisting this with $\Omega^1_Z \otimes L^{-1} \cong {\mathcal O}_Z(-3)$ gives an exact sequence
$$
0 \to {\mathcal O}_Z(-2) \to \Omega^1_{\overline X}(\log D) \otimes \Omega^1_Z \otimes L^{-1} \to {\mathcal O}_Z(-1) \to 0.
$$
On global sections this proves the assertion.
\end{proof}

So far we have only shown the vanishing of $H^q(X,{\mathbb V}_1)$ and hence of $H^q(X,{\mathbb V})$ for $q \ne 2$.
In order to apply the method of Gordon, Hanamura and Murre, we also have to deal with the case $\Lambda^i {\mathbb V}$.

\begin{theorem} \label{vanishing}
Let $\rho$ be an irreducible, non--constant representation of $\pi_1(X)$, which is a direct factor in
$\Lambda^k ({\mathbb V}_1 \oplus {\mathbb V}_2)$ for $k \le 2$. Then the intersection cohomology group
$$
H^1(X,{\mathbb V}_\rho)
$$
is zero.
\end{theorem}

\begin{proof} Let us first compute all such representations: if $k=1$ we have only ${\mathbb V}_1$ and its dual.
If $k=2$, we have the decomposition
$$
\Lambda^2 ({\mathbb V}_1 \oplus {\mathbb V}_2) = \Lambda^2 {\mathbb V}_1 \oplus
\Lambda^2 {\mathbb V}_2 \oplus {\rm End}({\mathbb V}_1).
$$
Since ${\mathbb V}_1$ is 3--dimensional, $\Lambda^2 {\mathbb V}_1  \cong {\mathbb V}_2$ and therefore
the only irreducible, non--constant representation that is new here
is ${\rm End}^0({\mathbb V}_1)$, the trace--free endomorphisms of ${\mathbb V}_1$.
Since we have already shown the vanishing $H^1(X,{\mathbb V}_{1,2})$, it remains to treat
$ H^1(X,{\rm End}^0({\mathbb V}_1))$.
The vanishing of $H^1(X,{\rm End}^0({\mathbb V}_1))$ is a general and well--known statement:
The representation ${\rm End}^0({\mathbb V}_1)$ has regular highest weight
and therefore contributes only to the middle dimension $H^2$. A reference for this is
\cite[Main Thm.]{li-sch}, cf. \cite[ch. VII]{borel-wallach} and \cite{vog-zuck}.
\end{proof}

The vanishing of $H^1(X,{\rm End}^0({\mathbb V}_1))$ has the following amazing consequence, which does not seem
easy to prove directly using purely algebraic methods. In the compact case this has been shown by Miyaoka, cf. \cite{miyaoka}.

\begin{lemma} In our situation we have
$$
H^0_{L^2}(\overline{X}, S^3 \Omega^1_{\overline{X}}(\log D)(-D) \otimes L^{-3})=0.
$$
\end{lemma}

\begin{proof} Write down the Higgs complex for ${\rm End}^0(E)$.
In degree one, a term which contains
$$
S^3 \Omega^1_{\overline{X}}(\log D)(-D) \otimes L^{-3}
$$
occurs. Since $H^1$ vanishes, this cohomology group must vanish too.
\end{proof}

Finally we want to discuss the case $k=3$. Unfortunately here the vanishing
techniques do not work in general. But we are able to at least give a bound for the dimension
of the remaining cohomology group. Namely for $k=3$, one has
$$
\Lambda^3 ({\mathbb V}_1 \oplus {\mathbb V}_2) = \Lambda^3 {\mathbb V}_1 \oplus \Lambda^3 {\mathbb V}_2
\oplus (\Lambda^2 {\mathbb V}_1 \otimes {\mathbb V}_2) \oplus (\Lambda^2 {\mathbb V}_2 \otimes {\mathbb V}_1).
$$
Here the only new irreducible and non--constant representation is
$$
S^2 {\mathbb V}_1 \subseteq  {\mathbb V}_1 \otimes {\mathbb V}_1
$$
and its dual. We would like to compute $H^1(X,S^2 {\mathbb V}_1)$
using a variant of the symmetric product of the $L^2$--complexes $\Omega^*(S)_{(2)}$ as
described in the appendix.
The Higgs complex without $L^2$--conditions looks as follows:
\begin{tiny}
$$
\begin{matrix}
\left(S^2 \Omega^1_{\overline X}(\log D) \otimes L^{-2} \right) \oplus \left( \Omega^1_{\overline X}(\log D) \otimes L^{-2} \right)
\oplus L^{-2} \cr \downarrow  \cr
\left( S^2 \Omega^1_{\overline X}(\log D) \otimes L^{-2} \otimes \Omega^1_{\overline X}(\log D) \right)
\oplus \left( \Omega^1_{\overline X}(\log D) \otimes L^{-2} \otimes \Omega^1_{\overline X}(\log D) \right)  \oplus
\left(L^{-2}  \otimes \Omega^1_{\overline X}(\log D) \right) \cr \downarrow  \cr
\left( S^2 \Omega^1_{\overline X}(\log D) \otimes L^{-2} \otimes \Omega^2_{\overline X}(\log D) \right)
\oplus \left( \Omega^1_{\overline X}(\log D) \otimes L^{-2} \otimes \Omega^2_{\overline X}(\log D) \right)  \oplus
\left(L^{-2}  \otimes \Omega^2_{\overline X}(\log D) \right)
\end{matrix}
$$
\end{tiny}
Again many pieces of differentials in this complex are isomorphisms or zero. For example the differential
$$
S^2 \Omega^1_{\overline X}(\log D) \otimes L^{-2} \otimes \Omega^1_{\overline X}(\log D)
\to \Omega^1_{\overline X}(\log D) \otimes L^{-2} \otimes \Omega^2_{\overline X}(\log D)
$$
is a projection map onto a direct summand, since for every vector space $W$ we have the
identity
$$
S^2 W \otimes W = S^3 W \oplus \left( W \otimes \Lambda^2 W \right).
$$
Therefore the Higgs complex for $S^2(E)$ is quasi--isomorphic to
$$
L^{-2} {\buildrel 0 \over \to}  S^3 \Omega^1_{\overline X}(\log D) \otimes L^{-2}
{\buildrel 0 \over \to} S^2 \Omega^1_{\overline X}(\log D) \otimes L^{-2}
\otimes \Omega^2_{\overline X}(\log D).
$$
We conclude that the first cohomology is given by
$$
H^0({\overline X},S^3 \Omega^1_{\overline X}(\log D) \otimes L^{-2}).
$$
If we additionally impose the $L^2$--conditions (see appendix), then we see that the first Higgs
cohomology of $S^2(E,\vartheta)$ vanishes, provided that we have
$$
H^0({\overline X},S^2 \Omega^1_{\overline X}(\log D) \otimes \Omega^1_{\overline X} \otimes L^{-2})=0.
$$
Using
$$
0 \to \sigma^* \Omega^1_{E \times E} \to \Omega^1_{\overline{X}} \to i_* \Omega^1_Z \to 0
$$
we obtain an exact sequence
$$
0 \to H^0(\overline{X}, S^2 \Omega^1_{\overline{X}}(\log D) \otimes L^{-2}) \to
H^0({\overline X},S^2 \Omega^1_{\overline X}(\log D) \otimes \Omega^1_{\overline X} \otimes L^{-2}) \to
$$
$$
\to H^0(Z,S^2 \Omega^1_{\overline X}(\log D) \otimes \Omega^1_Z \otimes L^{-2}).
$$
A generalization of \cite[example 3]{miyaoka} leads to the vanishing
$$
H^0(\overline{X}, S^2 \Omega^1_{\overline{X}}(\log D) \otimes L^{-2})=0.
$$
Since $\Omega^1_{\overline{X}}(\log D)|_Z = {\mathcal O}_Z(1) \oplus {\mathcal O}_Z(2)$, we get
$$
H^0(Z,S^2 \Omega^1_{\overline X}(\log D) \otimes \Omega^1_Z \otimes L^{-2})=\C^3,
$$
because $\Omega^1_{\overline X}(\log D) \otimes \Omega^1_Z \otimes L^{-2} = {\mathcal O}_Z \oplus
2{\mathcal O}_Z(-1) \oplus {\mathcal O}_Z(-2)$.
However we are not able to decide whether these $3$ sections lift to $\overline{X}$.
When we restrict to forms with fewer poles, then the vanishing will hold for a kind of
cuspidal cohomology.

\section{Proof of the Main Theorem}
\label{main}

In paragraph \ref{withmotdec} we prove our main theorem, in
paragraph \ref{goal} we give some indication on the proof of the
motivic decomposition conjecture in our case, however the details
will be published in a forthcoming paper. We thus will drop the
assumption on the motivic decomposition conjecture in Theorem
\ref{maintheorem}.

\subsection{From Relative to Absolute}
\label{withmotdec}

We now state and prove our main theorem. Let $p: {\overline \ca}
\longrightarrow {\overline X}$ be the compactified family over
Holzapfel's surface. Assume the motivic decomposition conjecture
\ref{MDC} (\cite{corti-hanamura}, \cite[Conj.
1.4]{gordon-hanamura-murre-i}) for ${\overline \ca}/{\overline X}$.
In the proof we will need an auxiliary statement which was
implicitely proven in section 6:

\begin{lemma} \label{fixpartlemma} Let $x \in X$ be a base point. Then $\pi_1^{\rm top}(X,x)$ acts on
the Betti cohomology group $H^{2j}(\mathcal
A_x(\mathbb{C}),\mathbb{Q})$. Then, for $0\le j\le d=3,$ the cycle
class map $CH^j(\mathcal{A})\to H^{2j}(\mathcal
A_x(\mathbb{C}),\mathbb{Q})^{\pi_1^{top}(X,x)}$ is surjective.
\end{lemma}

\begin{proof} By Lemma 5.3 the sheaf $R^1 p_*\C$ is a sum of two irreducible representation of 
$\pi_1^{top}(X,x)$. By the proof of Theorem 6.5., $R^2 p_*\C$ decomposes into a one--dimensional constant representation and three irreducible ones. The constant part corresponds to the identity in ${\rm End}(\V_1)$ 
and therefore to the polarization class on the fibers, which is a Hodge class.  Therefore the invariant classes
in $H^{2}(\mathcal A_x(\mathbb{C}),\mathbb{Q})$, and by duality also in $H^{4}(\mathcal A_x(\mathbb{C}),\mathbb{Q})$, consist of Hodge classes and are hence in the image of the cycle class map by the Hodge conjecture for divisors (and curves).
\end{proof}

Now we can prove our main theorem:

\begin{theorem}
\label{maintheorem} Assuming the motivic decomposition conjecture
\ref{MDC}, the total space of the family $p: {\overline \ca}
\longrightarrow {\overline X}$ supports a partial set of
Chow--K\"unneth projectors $\pi_i$ for $i \neq 4,5,6$.
\end{theorem}

\begin{proof} The motivic decomposition conjecture \ref{MDC} states that
we have a relative Chow--K\"unneth decomposition with projectors
$\Pi^i_\alpha$ on strata $X_\alpha$ which is compatible with the
topological decomposition theorem \cite{BBD}
$$
\sum_{j,\alpha} \Psi^j_\alpha: \R\overline{p}_*\Q_{\overline{A}}
{\buildrel \cong \over \to} \bigoplus_{j,\alpha} IC_{X_\alpha}(
{\mathcal V}^j_\alpha)[-j-\dim(X_\alpha)].
$$
Now we want to pass from relative Chow--K\"unneth decompositions to
absolute ones. We use the notation of
\cite{gordon-hanamura-murre-ii} and for the reader's convenience we
recall everything. Let $P^i/X$ and $P^i_\alpha/X$ be the mutually
orthogonal projectors adding up to the identity $\Delta(A/X) \in
CH_{\dim(A)}(A \times_X A)$ such that
$$
(P^i/X)_* \R p_*\Q_A = IC_X(R^i p_*\Q_A)[-i], \; (P^i_\alpha/X)_*\R
p_*\Q_A= IC_{X_\alpha}({\mathcal V}_\alpha^i)[-i-\dim(X_\alpha)],
$$
where the sheaves ${\mathcal V}_\alpha^i$ are local systems
supported over the cusps. The projectors $P^i_\alpha/X$ on the
boundary strata decompose further into Chow--K\"unneth components,
since the boundary strata consist of smooth elliptic curves and the
stratification has the product type fibers described in
theorem~\ref{tor.comp}. Let us now summarize what we know about the
local systems $R^ip_*\C$ on the open stratum $X_0$ from section 6:
$R^1 p_*\C$ is a sum of two irreducible representations 
and has no cohomology except in degree
$2$ by Theorem~\ref{vanishing}. $R^2 p_*\C$ contains a trivial
subsystem and the remaining complement has no cohomology except in
degree $2$ again by Theorem~\ref{vanishing}. $R^3 p_*\C$ also
contains a trivial subsystem and its complement has cohomology
possibly in degrees $1,2,3$, see section 6. By duality similar
properties hold for $R^ip_*\C$ with $i=4,5,6$. Using these
properties together with Lemma \ref{fixpartlemma} we can follow
closely the proof of Thm. 1.3 in \cite{gordon-hanamura-murre-ii}:
First construct projectors $ (P^{2r}/X)_{\rm alg}$ which are
constituents of $(P^{2r}/X)$ for $0 \le r \le 3$. This follows
directly from Lemma \ref{fixpartlemma} as in Step II of
\cite[section 1.7.]{gordon-hanamura-murre-ii}. Step III from
\cite[section 1.7.]{gordon-hanamura-murre-ii} is valid by the
vanishing observations above. As in Step IV of loc. cit. this
implies that we have a decomposition into motives in ${\rm
CH}{\mathcal M}(k)$
$$
M^{2r-1}=(A,P^{2r-1},0) \; (1 \le r \le d), \; M^{2r}_{\rm
trans}=(A,P^{2r}_{\rm trans},0) \; (0 \le r \le d),
$$
$$
M^{2r}_{\rm alg}=(A,P^{2r}_{\rm alg},0) \; (0 \le r \le d),
$$
plus additional boundary motives $M^j_\alpha$ for each stratum
$X_\alpha$. As in Step V of \cite{gordon-hanamura-murre-ii} we can
split $M^{2r}_{\rm alg}$ further. The projectors constructed in this
way define a set of Chow--K\"unneth projectors $\pi_i$ for $i \neq
4,5,6$, since the relative projectors which contribute to more than
one cohomology only affect cohomological degrees $4,5$ and $6$.
\end{proof}

\begin{rem}
If $H^1(X,S^2{\mathbb V}_1)$ vanishes or consists of algebraic
$(2,2)$ Hodge classes only, then we even obtain a complete
Chow--K\"unneth decomposition in the same way, since algebraic Hodge
$(p,p)$--classes define Lefschetz motives $\Z(-p)$ which can be
split off by projectors in a canonical way. Therefore the Hodge
conjecture on $A$ would imply a complete Chow--K\"unneth
decomposition. However the Hodge conjecture is not very far from
proving the total decomposition directly.
\end{rem}

\subsection{Motivic decomposition conjecture}
\label{goal}

The goal of this paragraph is to sketch the proof of the motivic
decomposition conjecture \ref{MDC} in the case we treat in this
paper. The complete details for the following argument will be
published in a future publication. First note that since $\ca$ is an
abelian variety we can use the work of Deninger and Murre
(\cite{den-mur}) on Chow-K\"unneth decompositions of Abelian schemes
to obtain relative Chow-K\"unneth projectors for $\ca /X$. To
actually get relative Chow-K\"unneth projectors for $\overline{\ca}
/\overline{X}$, we observe the following.

Recall our results in section~\ref{sec:compact}. We showed that the
special fibres over the smooth elliptic cusp curves $D_i$ are of the
form $Y_s=E \times \mathbb P^1\times\mathbb P^1$. We do not need the
cycle class map
$$
CH_*(Y_s \times Y_s) \to H_*(Y_s \times Y_s)
$$
to be an isomorphism as in \cite[Thm. I]{gordon-hanamura-murre-i}.
Since the boundary strata on ${\overline X}$ are smooth elliptic
curves it is sufficient to know the Hodge conjecture for the special
fibres. But the special fibres are composed of elliptic curves and
rational varieties by our results in section~\ref{sec:compact}.
Therefore the methods in \cite{gordon-hanamura-murre-i} can be
refined to work also in this case and we can drop the assumption in
theorem~\ref{maintheorem}.

\begin{remarks}
We hope to come back to this problem later and prove the motivic
decomposition conjecture for all Picard families. The existence of
absolute Chow--K\"unneth decompositions however seems to be out of
reach for other examples since vanishing results will hold only for
large arithmetic subgroups, i.e., small level.
\end{remarks}

\section{Appendix: Algebraic $L^2-$sub complexes of symmetric powers of
the uniformizing bundle of a two--dimensional complex ball quotient}
\label{appendix}
\newcommand{\stk}[1]{\stackrel{#1}{\longrightarrow}}
\newcommand{\IR}{\mbox{$I\!\!R$}}
\newcommand{\ZZ}{\mbox{$Z\!\!Z$}}
\newcommand{\IC}{\mbox{$C\!\!\!I$}}
\newcommand{\IH}{\mbox{$I\!\!H$}}
\newcommand{\IQ}{\mbox{$I\!\!Q$}}
\newcommand{\IO}{\mbox{$I\!\!O$}}
\newcommand{\IP}{\mbox{$I\!\!P$}}

$\overline X$ a $2$-dim projective variety with a normal crossing
divisor $D$, $X={\overline X}\setminus D$; assume that the
coordinates near the divisor are $z_1, z_2$.
\\

Consider the uniformizing bundle of a $2$-ball quotient
\[
E=\left( \Omega^1_{\overline X}(\log
D)\otimes\mathcal{K}^{-1/3}_{\overline X}(\log
D) \right) \oplus\mathcal{K}^{-1/3}_{\overline X}(\log D)
\]

We consider two cases: 1) $D$ is a smooth divisor (the case we need) and
2) $D$ is a normal crossing divisor.

\vskip .3cm {\bf Case 1}: Assume that $D$ is defined by $z_1=0$.
Taking $v$ as the generating section of
$\mathcal{K}^{-1/3}_{\overline X}(\log D)$,
${\frac{dz_1}{z_1}}\otimes v, dz_2\otimes v$ as the generating
sections of $\Omega^1_{\overline X}(\log
D)\otimes\mathcal{K}^{-1/3}_{\overline X}(\log D)$, then the Higgs
field
\[
\vartheta: E\to E \otimes\Omega^1_{\overline X}(\log D)
\]
is defined by setting $\vartheta({\frac{dz_1}{z_1}}\otimes
v)=v\otimes{\frac{dz_1}{z_1}}$, $\vartheta(dz_2\otimes v)=v\otimes
dz_2$, and $\vartheta(v)=0$.
\\

Clearly, if $\vartheta$ is written as
$N_1{\frac{dz_1}{z_1}}+N_2dz_2$, then
$N_1({\frac{dz_1}{z_1}}\otimes v)=v$, $N_1(dz_2\otimes v)=0$,
$N_1(v)=0$, $N_2({\frac{dz_1}{z_1}}\otimes v)=0$, $N_2(dz_2\otimes
v)=v$, $N_2(v)=0$; the kernel of $N_1$ is the subsheaf generated
by $dz_2\otimes v$ and $v$. Using the usual notation, we then have
\\
\begin{eqnarray*}
{\text{Gr}}_1W(N_1)&=& {\text {generated by}}~~
{\frac{dz_1}{z_1}}\otimes v\\
{\text{Gr}}_0W(N_1)&=& {\text {generated by}}~~
dz_2\otimes v\\
{\text{Gr}}_{-1}W(N_1)&=& {\text {generated by}} ~~v.
\end{eqnarray*}
So, one has
\begin{eqnarray*}
\Omega^0(E)_{(2)}&=& z_1\{{\frac{dz_1}{z_1}}\otimes
v\}+\{dz_2\otimes v\}+\{v\}\\
&=&{\text{Ker}}N_1+z_1 E;\\
\Omega^1(E)_{(2)}&=&{\frac{dz_1}{z_1}}\otimes(z_1\{{\frac{dz_1}{z_1}}\otimes
v\}+z_1\{dz_2\otimes v\}+z_1\{v\})\\
&+&dz_2\otimes(z_1\{{\frac{dz_1}{z_1}}\otimes
v\}+\{dz_2\otimes v\}+\{v\})\\
&=&{\frac{dz_1}{z_1}}\otimes z_1 E+dz_2\otimes({\text{Ker}}N_1+z_1 E);\\
\Omega^2(E)_{(2)}&=&{\frac{dz_1}{z_1}}\wedge
dz_2\otimes(z_1\{{\frac{dz_1}{z_1}}\otimes v\}+z_1\{dz_2\otimes
v\}+z_1\{v\})\\
&=&{\frac{dz_1}{z_1}}\wedge dz_2\otimes z_1 E,
\end{eqnarray*}
where $\{ x \}$ represents the line bundle generated by an element $x$.

\vskip .3cm {\bf Case 2}: As before, taking $v$ as the generating
section of $\mathcal{K}^{-1/3}_{\overline X}(\log D)$,
${\frac{dz_1}{z_1}}\otimes v, {\frac{dz_2}{z_2}}\otimes v$ as the
generating sections of $\Omega^1_{\overline X}(\log
D)\otimes\mathcal{K}^{-1/3}_{\overline X}(\log D)$, then the Higgs
field
\[
\vartheta: E \to E \otimes\Omega^1_{\overline X}(\log D)
\]
is defined by setting $\vartheta({\frac{dz_1}{z_1}}\otimes
v)=v\otimes{\frac{dz_1}{z_1}}$, $\vartheta({\frac{dz_2}{z_2}}\otimes
v)=v\otimes{\frac{dz_2}{z_2}}$, and $\vartheta(v)=0$.
\\

Clearly, if $\vartheta$ is written as
$N_1{\frac{dz_1}{z_1}}+N_2{\frac{dz_2}{z_2}}$, then
$N_1({\frac{dz_1}{z_1}}\otimes v)=v$,
$N_1({\frac{dz_2}{z_2}}\otimes v)=0$, $N_1(v)=0$,
$N_2({\frac{dz_1}{z_1}}\otimes v)=0$,
$N_2({\frac{dz_2}{z_2}}\otimes v)=v$, $N_2(v)=0$; the kernel of
$N_1$ (resp. $N_2$) is the subsheaf generated by
${\frac{dz_2}{z_2}}\otimes v$ (resp. ${\frac{dz_1}{z_1}}\otimes
v$) and $v$. We then have
\begin{eqnarray*}
{\text{Gr}}_1W(N_1)&=& {\text {generated by}}~~
{\frac{dz_1}{z_1}}\otimes v\\
{\text{Gr}}_0W(N_1)&=& {\text {generated by}}~~
{\frac{dz_2}{z_2}}\otimes v\\
{\text{Gr}}_{-1}W(N_1)&=& {\text {generated by}} ~~v\\
{\text{Gr}}_1W(N_2)&=& {\text {generated by}}~~
{\frac{dz_2}{z_2}}\otimes v\\
{\text{Gr}}_0W(N_2)&=& {\text {generated by}}~~
{\frac{dz_1}{z_1}}\otimes v\\
{\text{Gr}}_{-1}W(N_2)&=& {\text {generated by}} ~~v\\
{\text{Gr}}_1W(N_1+N_2)&=& {\text {generated by}}~~
({\frac{dz_1}{z_1}}+{\frac{dz_2}{z_2}})\otimes v\\
{\text{Gr}}_0W(N_1+N_2)&=& {\text {generated by}}~~
({\frac{dz_1}{z_1}}-{\frac{dz_2}{z_2}})\otimes v\\
{\text{Gr}}_{-1}W(N_1+N_2)&=& {\text {generated by}} ~~v\\
\end{eqnarray*}
So, one has
\begin{eqnarray*}
\Omega^0(E)_{(2)}&=&z_1\{{\frac{dz_1}{z_1}}\otimes
v\}+z_2\{{\frac{dz_2}{z_2}}\otimes v\}+\{v\}\\
&=&{\text{Ker}}N_1\cap{\text{Ker}}N_2+z_2{\text{Ker}}N_1+z_1{\text{Ker}}N_2;\\
\Omega^1(E)_{(2)}&=&{\frac{dz_1}{z_1}}\otimes(z_1\{{\frac{dz_1}{z_1}}\otimes
v\}+z_1z_2\{{\frac{dz_2}{z_2}}\otimes v\}+z_1\{v\})\\
&+&{\frac{dz_2}{z_2}}\otimes(z_2\{{\frac{dz_2}{z_2}}\otimes
v\}+z_1z_2\{{\frac{dz_1}{z_1}}\otimes v\}+
z_2\{v\})\\
&=&{\frac{dz_1}{z_1}}\otimes(z_1{\text{Ker}}N_2+z_1z_2{\text{Ker}}N_1)
+{\frac{dz_2}{z_2}}\otimes(z_2{\text{Ker}}N_1+z_1z_2{\text{Ker}}N_2);\\
\Omega^2(E)_{(2)}&=&{\frac{dz_1}{z_1}}\wedge{\frac{dz_2}{z_2}}\otimes
z_1z_2 E,
\end{eqnarray*}
\\

For the above two cases, it is to easy to check that
$\vartheta(\Omega^0(E)_{(2)})\subset\Omega^1(E)_{(2)}$
and $\vartheta(\Omega^1(E)_{(2)})\subset\Omega^2(E)_{(2)}$.
Thus, together $\vartheta\wedge\vartheta=0$, we have the
complex $(\{\Omega^i(E)_{(2)}\}_{i=0}^2, \vartheta)$
\[
0\to\Omega^0(E)_{(2)}\to\Omega^1(E)_{(2)}\to\Omega^2(E)_{(2)}\to 0
\]
with $\vartheta$ as the boundary operator. \\\\

Now we take  the $2^{\text{nd}}$-order symmetric power of
$(E, \vartheta)$, we obtain a new Higgs bundle ${\mathcal{S}}^2(E, \vartheta)$
(briefly, the Higgs field is still denoted by $\vartheta$)
as follows,
\begin{eqnarray*}
{\mathcal{S}}^2(E, \vartheta)
&=&{\mathcal{S}}^2\big(\Omega^1_{\overline X}(\log
D)\big)\otimes\mathcal{K}^{-2/3}_{\overline X}(\log
D)\oplus\\
&&\oplus\Omega^1_{\overline X}(\log
D)\otimes\mathcal{K}^{-2/3}_{\overline X}(\log
D)\oplus\mathcal{K}^{-2/3}_{\overline X}(\log D).
\end{eqnarray*}
The Higgs field $\vartheta$ maps ${\mathcal{S}}^2(E, \vartheta)$
into ${\mathcal{S}}^2(E, \vartheta)\otimes\Omega^1_{\overline
X}(\log D)$ and ${\mathcal{S}}^2(E,
\vartheta)\otimes\Omega^1_{\overline X}(\log D)$ into
${\mathcal{S}}^2(E, \vartheta)\otimes\Omega^2_{\overline
X}(\log D)$ so that one has a complex with the differentiation
$\vartheta$ as follows
\[
(*)~~~~~~0\to{\mathcal{S}}^2(E,
\vartheta)\to{\mathcal{S}}^2(E,
\vartheta)\otimes\Omega^1_{\overline X}(\log
D)\to{\mathcal{S}}^2(E, \vartheta)\otimes\Omega^2_{\overline
X}(\log D)\to 0;
\]
more precisely, one has
\begin{eqnarray*}
&&\vartheta\big({\mathcal{S}}^2\big(\Omega^1_{\overline X}(\log
D)\big)\otimes\mathcal{K}^{-2/3}_{\overline X}(\log
D)\big)\subset\big(\Omega^1_{\overline X}(\log
D)\otimes\mathcal{K}^{-2/3}_{\overline X}(\log
D)\big)\otimes\Omega^1_{\overline X}(\log D)\\
&&\vartheta\big(\Omega^1_{\overline X}(\log
D)\otimes\mathcal{K}^{-2/3}_{\overline X}(\log
D)\big)\subset\mathcal{K}^{-2/3}_{\overline X}(\log
D)\otimes\Omega^1_{\overline X}(\log
D)\\
&&\vartheta\big(\mathcal{K}^{-2/3}_{\overline X}(\log D)\big)=0\\
\end{eqnarray*}
and
\begin{eqnarray*}
&&\vartheta\big(\big({\mathcal{S}}^2\big(\Omega^1_{\overline X}(\log
D)\big)\otimes\mathcal{K}^{-2/3}_{\overline X}(\log
D)\big)\otimes\Omega^1_{\overline X}(\log
D)\big)\\
&&\subset\big(\Omega^1_{\overline X}(\log
D)\otimes\mathcal{K}^{-2/3}_{\overline X}(\log D)
\big)\otimes\Omega^2_{\overline X}(\log
D)\\
&&\vartheta\big(\big(\Omega^1_{\overline X}(\log
D)\otimes\mathcal{K}^{-2/3}_{\overline X}(\log D)
\big)\otimes\Omega^1_{\overline X}(\log
D)\big)\subset\mathcal{K}^{-2/3}_{\overline X}(\log
D)\otimes\Omega^2_{\overline X}(\log
D)\\
&&\vartheta\big(\mathcal{K}^{-2/3}_{\overline X}(\log
D)\otimes\Omega^1_{\overline X}(\log
D)\big)=0\\
\end{eqnarray*}

\vskip .5cm \noindent Note: Let $V$ be a $SL(2)$-module, then
\[
\mathcal{S}^2V\otimes V\simeq\mathcal{S}^3V\oplus
V\otimes\wedge^2V.
\]
In general, one needs to consider the representations of $GL(2)$;
in such a case, we can take the determinant of the representation
in question, and then go back to a representation of $SL(2)$.

\begin{eqnarray*}
{\mathcal{S}}^2(E, \vartheta)
&=&{\mathcal{S}}^2\big(\Omega^1_{\overline X}(\log
D)\big)\otimes\mathcal{K}^{-2/3}_{\overline X}(\log
D)\\
&&\oplus\Omega^1_{\overline X}(\log
D)\otimes\mathcal{K}^{-2/3}_{\overline X}(\log
D)\\
&&\oplus\mathcal{K}^{-2/3}_{\overline X}(\log D)\\
{\mathcal{S}}^2(E, \vartheta)\otimes\Omega^1_{\overline X}(\log
D)&=&\big({\mathcal{S}}^2\big(\Omega^1_{\overline X}(\log
D)\big)\otimes\mathcal{K}^{-2/3}_{\overline X}(\log
D)\big)\otimes\Omega^1_{\overline X}(\log D)\\
&&\oplus\big(\Omega^1_{\overline X}(\log
D)\otimes\mathcal{K}^{-2/3}_{\overline X}(\log
D)\big)\otimes\Omega^1_{\overline X}(\log
D)\\
&&\oplus\mathcal{K}^{-2/3}_{\overline X}(\log
D)\otimes\Omega^1_{\overline X}(\log D)\\
 {\mathcal{S}}^2(E,
\vartheta)\otimes\Omega^1_{\overline X}(\log
D)&=&\big({\mathcal{S}}^2\big(\Omega^1_{\overline X}(\log
D)\big)\otimes\mathcal{K}^{-2/3}_{\overline X}(\log
D)\big)\otimes\Omega^2_{\overline X}(\log D)\\
&&\oplus\big(\Omega^1_{\overline X}(\log
D)\otimes\mathcal{K}^{-2/3}_{\overline X}(\log
D)\big)\otimes\Omega^2_{\overline X}(\log D)\\
&&\oplus\mathcal{K}^{-2/3}_{\overline X}(\log
D)\otimes\Omega^2_{\overline X}(\log D)
\end{eqnarray*}

Assuming that the divisor $D$ is smooth, we next want to consider
the $L^2$-holomorphic Dolbeault sub-complex of the above complex
(*):
\[
0\to({\mathcal{S}}^2(E,
\vartheta))_{(2)}\to({\mathcal{S}}^2(E,
\vartheta)\otimes\Omega^1_{\overline X}(\log
D))_{(2)}\to({\mathcal{S}}^2(E,
\vartheta)\otimes\Omega^2_{\overline X}(\log D))_{(2)}\to 0,
\]
and explicitly write down $({\mathcal{S}}^2(E,
\vartheta)\otimes\Omega^i_{\overline X}(\log D))_{(2)}$.
\\
Note that taking symmetric power for $L^2$-complex does not have
obvious functorial properties in general.

\vskip .5cm \noindent We will continue to use the previous
notations. For simplicity, we will further set
$v_1={\frac{dz_1}{z_1}}\otimes v$ and $v_2={dz_2}\otimes v$; we
also denote $e_1\otimes e_2+e_2\otimes e_1$ by $e_1\odot e_2$, the
symmetric product of the vectors $e_1$ and $e_2$.
\\
Thus, ${\mathcal{S}}^2\big(\Omega^1_{\overline X}(\log
D)\big)\otimes\mathcal{K}^{-2/3}_{\overline X}(\log D)$, as a
sheaf, is generated by $v_1\odot v_1$, $v_1\odot v_2$, $v_2\odot
v_2$; $\Omega^1_{\overline X}(\log
D)\otimes\mathcal{K}^{-2/3}_{\overline X}(\log D)$ is generated by
$v_1\odot v$, $v_2\odot v$; and $\mathcal{K}^{-2/3}_{\overline
X}(\log D)$ is generated by $v\odot v$. Also, it is easy to check
how $N_1, N_2$ act on these generators; as for $N_1$, we have
(Note $N_1v_1=v, N_1v_2=0, N_1v=0$.)
\begin{eqnarray*}
N_1(v_1\odot v_1)&=&2v_1\odot v\\
N_1(v_1\odot v_2)&=&v_2\odot v\\
N_1(v_2\odot v_2)&=&0\\
N_1(v_1\odot v)&=&v\odot v\\
N_1(v_2\odot v)&=&0\\
N_1(v\odot v)&=&0.
\end{eqnarray*}
Clearly, $N_1$ maps ${\mathcal{S}}^2\big(\Omega^1_{\overline
X}(\log D)\big)\otimes\mathcal{K}^{-2/3}_{\overline X}(\log D)$
into $\Omega^1_{\overline X}(\log
D)\otimes\mathcal{K}^{-2/3}_{\overline X}(\log D)$,
$\Omega^1_{\overline X}(\log
D)\otimes\mathcal{K}^{-2/3}_{\overline X}(\log D)$ into
$\mathcal{K}^{-2/3}_{\overline X}(\log D)$, and then
$\mathcal{K}^{-2/3}_{\overline X}(\log D)$ to $0$. So, $N_1$ is of
index $3$ on the $2^{\text{nd}}$-order symmetric power
$\mathcal{S}^2 E$(as is obvious from the abstract theory
since $N_1$ is of index $2$ on $E$); and we then have the
following gradings
\begin{eqnarray*}
{\text{Gr}}_2W(N_1)&=& {\text {generated by}}~~ v_1\odot v_1\\
{\text{Gr}}_1W(N_1)&=& {\text {generated by}}~~ v_1\odot v_2\\
{\text{Gr}}_0W(N_1)&=& {\text {generated by}}~~ v_1\odot v, v_2\odot v_2\\
{\text{Gr}}_{-1}W(N_1)&=& {\text {generated by}} ~~v_2\odot v\\
{\text{Gr}}_{-2}W(N_1)&=& {\text {generated by}}~~ v\odot v.\\
\end{eqnarray*}
(Note that $N_1$, acting on $E$, has two invariant
(irreducible) components, one being generated by $v_1, v$, the
other by $v_2$, so that $N_1$ has three invariant components on
$\mathcal{S}^2E$, as is explicitly showed in the above
gradings.)

\vskip .5cm\noindent Now we can write down $L^2$-holomorphic
sections of $\mathcal{S}^2 E$, namely the sections generated
by $v_1\odot v, v_2\odot v_2, v_2\odot v, v\odot v$, and
$z_1\mathcal{S}^2 E$; in the invariant terms, they should be
\[
(\mathcal{S}^2(E, \vartheta))_{(2)}=E\odot{\text{Im}}N_1+
\mathcal{S}^2({\text{Ker}}N_1)+z_1\mathcal{S}^2 E.
\]

\vskip .5cm\noindent Now it is easy to also write down
$({\mathcal{S}}^2(E, \vartheta)\otimes\Omega^1_{\overline
X}(\log D))_{(2)}$ and $({\mathcal{S}}^2(E,
\vartheta)\otimes\Omega^2_{\overline X}(\log D))_{(2)}$:
\begin{eqnarray*}
({\mathcal{S}}^2(E, \vartheta)\otimes\Omega^1_{\overline
X}(\log
D))_{(2)}&=&{\frac{dz_1}{z_1}}\otimes(\mathcal{S}^2({\text{Im}}N_1)+z_1\mathcal{S}^2 E)\\
&&+dz_2\otimes(E \odot{\text{Im}}N_1+\mathcal{S}^2({\text{Ker}}N_1)+z_1\mathcal{S}^2 E);\\
({\mathcal{S}}^2(E, \vartheta)\otimes\Omega^2_{\overline
X}(\log D))_{(2)}&=&{\frac{dz_1}{z_1}}\wedge
dz_2\otimes(\mathcal{S}^2({\text{Im}}N_1)+z_1\mathcal{S}^2 E).
\end{eqnarray*}

Similary, one can determine the algebraic $L^2-$sub complex of $S^n(E,\vartheta)$ for any $n\in\mathbb N.$

\section{Acknowledgement}

It is a pleasure to dedicate this work to Jaap Murre who has been so tremendously important for
the mathematical community. In particular we want to thank him for his constant support
during so many years.  \\
We are grateful to Bas Edixhoven, Jan Nagel and Chris Peters for organizing such a
wonderful meeting in Leiden. Many thanks go to F. Grunewald, R.-P. Holzapfel,
J.-S. Li and J. Schwermer for very helpful discussions.

\end{document}